\documentclass[10pt]{article}

\usepackage{amsthm,amsmath,amssymb}
\usepackage[isolatin]{inputenc}
\usepackage{lmodern}
\usepackage[T1]{fontenc}
\usepackage{url}

\newcommand{\xf}{\mathbf{x}}
\newcommand{\af}{\mathbf{a}}
\newcommand{\fs}{\mathbf{s}}
\newcommand{\vf}{\mathbf{v}}
\newcommand{\grad}{\mathop{\mathrm{grad}}}
\newcommand{\vt}{\tilde{V}}
\newcommand{\R}{\mathbb{R}}
\newcommand{\ed}{\mathrm{d}}
\newcommand{\vektor}[2]{\left[\begin{array}{c} #1 \\ #2
    \end{array}\right]} 
\newcommand{\vektort}[3]{\left[\begin{array}{c} #1 \\ #2 \\ #3
    \end{array}\right]}   
\newcommand{\matris}[4]{\left[\begin{array}{cc} #1 & #2 \\ #3 &
      #4\end{array}\right]} 
\newcommand{\matrist}[9]{\left[\begin{array}{ccc} #1 & #2 & #3 \\ #4 &
      #5 & #6 \\ #7 & #8 & #9\end{array}\right]}
\newcommand{\pd}[2]{\frac{\partial #1}{\partial #2}}

\newtheorem{Thm3}{Theorem}
\newtheorem{Lem3}[Thm3]{Lemma}
\newtheorem{Prop3}[Thm3]{Proposition}
\newtheorem{Rem3}{Remark}

\newtheorem{Exmp3}{Example}
\newtheorem{Hyp3}{Hypothesis}
\begin{document}
\bibliographystyle{abbrv}

\title{Multiplication for solutions of the equation $\grad{f} =
  M\grad{g}$}    

\author{{\bf Jens Jonasson}\\ Department of Mathematics\\ Link\"oping
  University\\ SE-581 83 Link\"oping, Sweden}   
\maketitle
\begin{center}
  E-mail address: \texttt{jejon@mai.liu.se} 
\end{center}

\begin{abstract}
Linear first order systems of partial differential equations of the
form $\nabla f = M\nabla g,$ where $M$ is a constant matrix, are
studied on vector spaces over the fields of real and complex numbers, respectively.
The Cauchy--Riemann equations belong to this class. We introduce a
bilinear $*$-multi\-plication on the solution space, which plays the
role of a nonlinear superposition principle, that allows for
algebraic construction of new solutions from known solutions. 
The gradient equations $\nabla f = M\nabla g$ constitute only a simple special case of a much larger class of systems of partial differential equations which admit a bilinear multiplication on the solution space, but we prove that any gradient equation has the exceptional property that the general analytic solution can be expressed through power
series of certain simple solutions, with respect to the
$*$-multiplication.  
\end{abstract}

\clearpage

\section{Introduction} 

We consider equations of the form 
\begin{align}\label{fmg}
  \nabla f = M\nabla g,
\end{align}
where $f(\xf)$ and $g(\xf)$ are unknown scalar-valued functions,
defined on an open convex domain in an $n$-dimensional vector space
$\mathcal{V}$ over the field of real or complex numbers, and
$M$ is a constant $n\times n$ matrix. Hence, the systems under
consideration constitute a system, that is usually overdetermined, of
$n$ first order linear partial differential equations (abbreviated
PDEs) for two unknown functions.

The Cauchy--Riemann equations
\begin{align*}
  \pd{f}{x} & = \pd{g}{y}\\
  \pd{f}{y} & = -\pd{g}{x}
\end{align*}
are an example of an equation of the form (\ref{fmg}). Since there is
a $1-1$ correspondence between solutions $(f,g)$ of the
Cauchy--Riemann equations and holomorphic functions $F = f + \mathrm{i}
g$, any product 
\begin{align*}
  (f + \mathrm{i} g)(\tilde{f} + \mathrm{i}\tilde{g}) = (f\tilde{f} -
  g\tilde{g}) + \mathrm{i}(f\tilde{g} + g\tilde{f})
\end{align*}
of two holomorphic functions is again holomorphic. The ordinary
multiplication of holomorphic functions defines a bilinear
$*$-multiplication on the solutions space $S$ of the Cauchy--Riemann
equations
\begin{align}\label{crmult}
  \begin{aligned}
    *: S \times S & \longrightarrow S \\
    (f, g)*(\tilde{f}, \tilde{g}) & = (f\tilde{f} - g\tilde{g},
    f\tilde{g} + g\tilde{f}).  
\end{aligned}
\end{align}
Moreover, since every holomorphic function is analytic, every solution
of the Cauchy--Riemann equations can be expressed locally as a power
series of a simple solution. For example, every solution can be formulated, in a neighborhood of the
origin, as a power series of the simple solution $(x, y)$
\begin{align}\label{crps}
  (f, g) = \sum_{r=0}^\infty (a_r, b_r)*(x, y)_*^r, \quad \textrm{where} \quad
  (x, y)_*^r = \underbrace{(x, y)* \cdots *(x, y)}_{r \textrm{ factors}},
\end{align}
and $a_r, b_r$ are real constants.

In \cite{jodeit-1990}, Jodeit and Olver have given an expression for
the general solution of the equation (\ref{fmg}). In
\cite{jonasson-2007}, we have have introduced a multiplication $*$ of
solutions for a wide class of linear first order systems of
differential equations containing systems of the form (\ref{fmg}),
where the matrix $M$ can also have non-constant entries. The
$*$-multiplication is a nonlinear superposition principle that allows
for algebraic construction of new solutions from known solutions. By
combining the $*$-multiplication with the ordinary linear
superposition principle, power series solutions can be constructed
from a single simple solution.

In this paper we show that the equations (\ref{fmg}) are included in
the family of systems equipped with a $*$-multiplication on the
solution space. We also compare power series solutions, constructed
with the $*$-multiplication, with the formulas in \cite{jodeit-1990}
that describe the general solution of (\ref{fmg}). This paper contains
the following main results:  
\begin{enumerate}
  \item Any equation (\ref{fmg}), including both the real and the
    complex case, admits a $*$-multiplication which provides a
    nonlinear superposition formula. For the Cauchy--Riemann
    equations, this multiplication reduces to the ordinary
    multiplication (\ref{crmult}), obtained from the multiplication of
    holomorphic functions.
  \item Any solution of (\ref{fmg}) can be represented as a power
    series (with respect to the $*$-multiplication) of certain simple
    solutions. In other words, we provide a different and more
    explicit way to describe the general analytic solution of
    (\ref{fmg}). This may be compared with the Cauchy--Riemann
    equations where each solution can be obtained from the
    harmonically conjugated pair given by the real and imaginary parts 
    of a power series in one complex variable (\ref{crps}).    
  \item The equations of the form (\ref{fmg}) constitute an
    interesting example of a large class of systems of linear PDEs with a
    simple structure, admitting $*$-multiplication of solutions. The
    content of this paper can be considered a thorough investigation of
    the $*$-multiplication for the class of systems of the form
    (\ref{fmg}). The results contribute to a better understanding of the
    general class of equations admitting $*$-multiplication. 
\end{enumerate}

\begin{Rem3}
  In the special case when the matrix $M$ consists of only one Jordan
  block for each eigenvalue, the general solution of $($\ref{fmg}$)$
  can be described in terms of component functions of functions which
  are differentiable over some algebra \cite{waterhouse-1992}. For
  this restricted class of matrices $M$, some of the results above can
  also be derived from the content of \cite{waterhouse-1992}.     
\end{Rem3}

This paper is organized as follows. In section \ref{multsec} we
present the necessary results from \cite{jonasson-2007} about
bilinear $*$-multiplication for solutions of linear systems of PDEs
with variable and constant coefficients. Section \ref{fmgsec} contains
a summary of the main results of paper \cite{jodeit-1990} by Jodeit
and Olver about the algebraic form of the general solution to the
equation (\ref{fmg}). Sections \ref{cstmultsec} and \ref{pssec}
contain the main results of this paper. In section \ref{cstmultsec},
we show that any equation of the form (\ref{fmg}) admits a bilinear
$*$-multiplication of solutions of the type described in
\cite{jonasson-2007}. In section \ref{pssec} we prove (theorem
\ref{ps2} and theorem \ref{realps}) that every solution of (\ref{fmg})
can be expressed as a power series of simple solutions with respect to
the $*$-multiplication.

\section{Multiplication of solutions for systems of
  PDEs}\label{multsec} 

Certain systems of linear partial differential equations admit a
bilinear operation on the solution space. The most general results
about this operation, which we call $*$-multiplication, are given in
\cite{jonasson-2007} (the results are derived for equations on a real
differentiable manifold, but all arguments hold also for analytic
equations over a complex vector space), and in this section we give a
brief summary. 

Let $Z_\mu = Z_0 + Z_1\mu + \cdots + Z_{m-1} \mu^{m-1} + \mu^m$ and
$V_\mu = V_0 + V_1\mu + \cdots + V_{m-1} \mu^{m-1}$ be polynomials in
the variable $\mu$, with coefficients that are smooth scalar-valued
functions on the vector space $\mathcal{V}$. Furthermore, let $A_\mu =
A_0 + A_1\mu + \cdots + A_k\mu^k$ be a polynomial where each
coefficient is a $n \times n$ matrix, with smooth scalar-valued
functions as entries. Consider then the matrix equation
\begin{align}\label{avz}
  A_\mu \nabla V_\mu \equiv 0 \quad \left(mod\; Z_\mu \right), 
\end{align}
where the unknown function $V_\mu$ is a solution if the remainder of
the polynomial $A_\mu \nabla V_\mu$ modulo $Z_\mu$ is zero, i.e., if
there exists a polynomial vector $\mathbf{u}_\mu$ such that $A_\mu
\nabla V_\mu = Z_\mu \mathbf{u}_\mu$.

Given two solutions $V_\mu$ and $W_\mu$, the $*$-product $V_\mu *
W_\mu$ is defined as the unique remainder of the ordinary product
$V_\mu W_\mu$ modulo $Z_\mu$, i.e., $V_\mu * W_\mu$ is the unique
polynomial of degree less than $m$ that can be written as $V_\mu *
W_\mu = V_\mu W_\mu - Q_\mu Z_\mu$ for some polynomial $Q_\mu$. In
general, neither $V_\mu W_\mu$ nor $V_\mu * W_\mu$ are solutions of
(\ref{avz}), but when $A_\mu$ and $Z_\mu$ are related by the equation
\begin{align}\label{azz}
  A_\mu \nabla Z_\mu \equiv 0 \quad \left(mod\; Z_\mu \right), 
\end{align}
that is, when $Z_\mu - \mu^m$ is a solution of (\ref{avz}), then the
$*$-multiplication maps any two solutions into a new solution. Thus,
the $*$-multiplication provides a method for constructing, in a pure
algebraic way, new solutions from already known solutions. Especially,
when the coefficients of $Z_\mu$ are not all constant, $*$-products of
trivial (constant) solutions are in general non-trivial solutions. 
   
\begin{Exmp3}
  If we let
  \begin{align*}
    A_\mu = \matris{\mu}{1}{-1}{\mu}, \quad V_\mu = g + f\mu, \quad
    Z_\mu = 1 + \mu^2,
  \end{align*}
  the equation $($\ref{avz}$)$ reduces to the Cauchy--Riemann
  equations. Since the coefficients in $Z_\mu$ are here constant, the
  condition $($\ref{azz}$)$ is trivially satisfied. Hence, the
  $*$-multiplication maps solutions to solutions, and it reconstructs
  the multiplication formula $($\ref{crmult}$)$ obtained from the
  multiplication of holomorphic functions:
  \begin{align*}
  (f + g\mu)(\tilde{f} + \tilde{g}\mu) & = f\tilde{f} + (f\tilde{g} +
  g\tilde{f})\mu + g\tilde{g}\mu^2 \\
  & \equiv (f\tilde{f} - g\tilde{g}) + (f\tilde{g} + g\tilde{f})\mu =
  (f + g\mu)*(\tilde{f} + \tilde{g}\mu).     
  \end{align*}
\end{Exmp3}

We see that, in order to obtain the first of the main results named in
the introduction, it is enough to find, for each matrix $M$, a matrix
$A_\mu$ and a function $Z_\mu$ such that the relation (\ref{azz}) is
satisfied, and such that the equation (\ref{avz}) is equivalent to
(\ref{fmg}). 

The theory of $*$-multiplication can also be formulated in matrix
language, without use of congruences and the extra parameter $\mu$. In
matrix notation the equation (\ref{avz}) can be written as 
\begin{align}\label{matrissys}
  \sum_{i=0}^k C^iV'A_i = 0,
\end{align}
where $V'$ denotes the $m \times n$ functional matrix with entries
$(V')_{ij} = \partial V_i / \partial x_j$ and $C$ is the companion
matrix for the polynomial $Z_\mu$ 
\begin{align*}
  C = \left[
  \begin{array}{ccccc}
    0 & 0 & \cdots & 0 & -Z_0 \\
    1 & 0 & \cdots & 0 & -Z_1 \\
    0 & 1 &  \ddots & & \vdots \\
    \vdots & \vdots & \ddots & 0 & -Z_{n-2} \\
    0 & 0 & \cdots & 1 & -Z_{n-1} 
  \end{array}
  \right].
\end{align*}
The $*$-multiplication exists for any equation (\ref{matrissys}) for
which $V = [Z_0 \quad Z_1 \;\; \cdots Z_{m-1}]^T$ is a solution.
The $*$-product of two solutions can then be written explicitly as the
following matrix product
\begin{align*}
  V*W & = V_CW_Ce_1 = \left( \sum_{i=0}^{m-1} V_iC^i \right)\left(
    \sum_{i=0}^{m-1} W_iC^i \right)\left[ \begin{array}{c} 1 \\ 0 \\
      \vdots \\ 0 \end{array} \right]. 
\end{align*}

An interesting property of the $*$-multiplication is that complex
solutions can be generated, starting with simple solutions.
Especially, when $Z_\mu$ has non-constant coefficients, $*$-powers of
trivial (constant) solutions will in general be non-trivial. By using
the $*$-multiplication in combination with the linear superposition
principle, one can form $*$-polynomial solutions
\begin{align*}
  V_\mu = \sum_{r=0}^N a_r\mu^r_*, \quad {\rm where} \quad \mu_*^r =
  \underbrace{\mu * \mu * \cdots *\mu,}_{r factors}   
\end{align*}
or in matrix notation
\begin{align}\label{polysol}
  V = \sum_{r=0}^N a_rC^re_1, 
\end{align}
where $a_i$ are constants. Moreover, when some conditions are
satisfied, by letting $N\rightarrow \infty$ in (\ref{polysol}) one
obtains power series solutions. However, when the coefficients of
$Z_\mu$ are all constant (which is the case for the Cauchy--Riemann
equations), power series of trivial solutions will be trivial.
Instead, in order to construct interesting solutions, one has to
build power series from a non-trivial solution. Such power series
are not considered in \cite{jonasson-2007} but the results for power
series of trivial solutions extend immediately to power series of
non-trivial solutions. More precisely, the power series $\sum
a_r(V_\mu)_*^r$, which in matrix notation reads
\begin{align*}
  \sum_{r=0}^\infty a_r\tilde{C}^re_1, \quad  {\rm where} \quad
  \tilde{C} = \sum_{i=0}^{m-1}V_iC^i,  
\end{align*}
defines a solution of (\ref{avz}) as long as the modulus of the
eigenvalues of $\tilde{C}$ are smaller than the radius of convergence
of the power series $\sum a_rt^r$ and the geometrical multiplicities
of the eigenvalues are constant. For power series of trivial
solutions, the matrix $\tilde{C}$ is identical with the companion
matrix $C$ and its eigenvalues coincide with the roots of the
polynomial $Z_\mu$. 

We can also consider $*$-power series with arbitrary trivial solutions
as coefficients 
\begin{align*}
  \sum_{r=0}^\infty a_{\mu, r}*(V_\mu)_*^r, \quad  {\rm where} \quad
  a_{\mu, r} = a_{r,0} + a_{r,1}\mu + \cdots a_{r,m-1}\mu^{m-1}.
\end{align*}

\begin{Rem3}
The $*$-multiplication \cite{jonasson-2006, jonasson-2007} is a
generalization of the multiplication of cofactor pair systems
introduced by Lundmark \cite{lundmark-2001}, which in turn is a
generalization of a recursive formula for obtaining new cofactor pair
systems \cite{rauch-1986, rauch-1999, lundmark-2003}. A cofactor pair
system is a Newton equation $\ddot{q}^h +
\Gamma^h_{ij}\dot{q}^i\dot{q}^j = F^h,$ $h = 1, 2, \ldots, n$ on a
Riemannian manifold, where the vector field $F$ can be written as
\begin{align*}
F = -(\det{J})^{-1}J\nabla V =
-(\det{\tilde{J}})^{-1}\tilde{J}\nabla\tilde{V} 
\end{align*}
for some functions $V,$ $\tilde{V}$ and two special conformal Killing
tensors $J$, $\tilde{J}$. The family of cofactor pair systems contain
all separable conservative Lagrangian systems, that in general can be
integrated through separation of variables in the Hamilton--Jacobi
sense \cite{rauch-1999,lundmark-2003,lundmark-2002,rauch-2003,
crampin-2001, benenti-2005}. The multiplication of cofactor pair
systems is a mapping that, for fixed special conformal Killing tensors
$J$ and $\tilde{J}$, for two solutions $(V, \tilde{V})$ and $(W,
\tilde{W})$ of the equation $(\det{J})^{-1}J\nabla V =
(\det{\tilde{J}})^{-1}\tilde{J}\nabla\tilde{V}$ prescribes a new
solution $(V, \tilde{V})*(W, \tilde{W})$ in a bilinear way.
\end{Rem3}

\section{The general solution of  $\nabla f = M\nabla g$}\label{fmgsec}

In \cite{jodeit-1990}, Jodeit and Olver describe a general analytic
solution of the equation (\ref{fmg}) for any constant, real or
complex, matrix $M$. In this section, we introduce some notation, and
briefly describe the main results from \cite{jodeit-1990}.
  
A linear change of variables $\xf\rightarrow A\xf$ transforms the
equation (\ref{fmg}) into $\nabla f = A^{-T}MA^T\nabla g$. Thus, by
performing this kind of transformation, we can assume that the matrix
$M$ is in some suitable canonical form with respect to similarity. We
will assume in the following that $M$ is in Jordan canonical form when
$\mathcal{V}$ is a vector space over the complex numbers, and in real
Jordan canonical form \cite{horn} for real $\mathcal{V}$.

\subsection{The complex case} 

If we let $\lambda_1,\ldots ,\lambda_p$ denote the distinct
eigenvalues of the matrix $M$, there is a primary decomposition of the
vector space  
\begin{align}\label{v1vp}
  \mathcal{V} = \mathcal{V}^1\oplus\cdots\oplus \mathcal{V}^p 
\end{align}
consisting of invariant subspaces $\mathcal{V}^k = \ker{(M -
\lambda_kI)^{e_k}}$, where $e_k$ is the geometric multiplicity of
$\lambda_k$. As mentioned above, there is no loss of generality to
assume that $M$ has a diagonal block structure $M =
\mathrm{diag}(M^1,\ldots ,M^p),$ where $M^k =
\mathrm{diag}(J^k_1,\ldots ,J^k_{p_k})$ and each $J^k_i$ is a Jordan
block corresponding to the eigenvalue $\lambda_k$, i.e., 
\begin{align*}
  J^k_i = \left[ 
    \begin{array}{cccc}
      \lambda_k & 1 & &\\
      & \lambda_k & \ddots &\\
      & & \ddots & 1\\
      & & & \lambda_k
    \end{array} \right].
\end{align*} 
There is a simple decomposition of the general solution of (\ref{fmg})
that is associated with the primary decomposition (\ref{v1vp}) of the
vector space $\mathcal{V}$. Every solution $(f,\; g)$ of the equation
(\ref{fmg}) can be written as a sum
\begin{align}\label{sumfg}
  f = f^1 + f^2 + \cdots + f^p, \quad g = g^1 + g^2 + \cdots + g^p,
\end{align}
where $(f^k,\; g^k)$ is a solution of the corresponding equation
$\nabla f = M^k\nabla g$. Hence, the problem of describing the general
solution of (\ref{fmg}) is reduced to the case when the matrix $M$
consists of Jordan blocks corresponding to a single eigenvalue
$\lambda$. Therefore, let $M = \mathrm{diag}(J_1,\ldots ,J_m),$ where
$J_k$ is a Jordan block of size $(n_k+1)\times (n_k+1)$ corresponding
to $\lambda$. We assume also that coordinates are chosen in such a way
that the size of the Jordan blocks is decreasing, i.e., $n_1 \ge n_2
\ge \cdots \ge n_m \ge 0$. Any vector $\xf\in \mathcal{V}$ is
decomposed as $\xf = [\xf^1 \quad \xf^2 \;\; \cdots \;\; \xf^m]^T$ so that
$\xf^k = [x^k_0 \quad x^k_1 \;\; \cdots \;\; x^k_{n_k}]^T$ contains the
variables that correspond to the block $J_k$. The variables $x^1_0,\;
x^2_0, \ldots, x^m_0$, that appear at the first position in the
vectors $\xf^1,\; \xf^2, \ldots, \xf^m$, are called \emph{major
variables} and all other variables are referred to as \emph{minor
variables}. For each block $J_k$, we introduce a function
\begin{align*}
  x^k_\mu = x_0^k + x_1^k\mu + \cdots + x_{n_k}^k\mu^{n_k},
\end{align*}
which depends on the variables contained in $\xf^k$ and a real
parameter $\mu$. Furthermore, let $\nu_k$ denote the number of Jordan
blocks of $M$ which are of size at least $k+1$, i.e., $\nu_k =
\mathrm{max}\{i\in\mathbb{N}:n_i \ge k\}$, and define vectors
$\xf^{(0)}_\mu, \ldots, \xf^{(n_1)}_\mu$ in the
following way:
\begin{align*}
  \mathbf{\xf}^{(k)}_\mu = [x^1_\mu \quad x^2_\mu \;\; \cdots \;\;
  x^{\nu_k}_\mu]^T.  
\end{align*}  
According to \cite{jodeit-1990}, two analytic functions $f(\xf)$ and
$g(\xf)$ form a solution of the equation (\ref{fmg}) (where we have
assumed that the matrix $M$ has only one eigenvalue $\lambda$ and has
a Jordan canonical form) if and only if
\begin{align*}
  f = f_1 + f_2 + \cdots + f_{n_1} + c, \quad g = g_1 + g_2 + \cdots +
  g_{n_1},  
\end{align*}  
where $c$ is an arbitrary constant, and there exist scalar-valued
analytic functions  
\begin{align*}
  \phi_k(s_1, \ldots, s_{\nu_k}),\quad k = 0, 1, \ldots, n_1,
\end{align*} 
such that $f_k$, $g_k$ are given by
\begin{align}\label{gensol}
  \left\{
    \begin{aligned}
      f_k(\xf) = & \left. \lambda\frac{\partial^k}{\partial
          \mu^k}\phi_k \left( \xf^{(k)}_\mu\right) \right|_{\mu=0} + 
      k\left.\frac{\partial^{k-1}}{\partial \mu^{k-1}}\phi_k \left(
          \xf^{(k)}_\mu\right) \right|_{\mu=0} \\ 
      g_k(\xf) = & \left. \frac{\partial^k}{\partial \mu^k}\phi_k
        \left( \xf^{(k)}_\mu\right) \right|_{\mu=0}. 
    \end{aligned}
  \right.
\end{align}   
There are two interesting observations regarding the general solution:  
\begin{enumerate}
  \item By a change of dependent variable $h = f - \lambda g$
    the equation (\ref{fmg}) transforms into an equation for
    $h$ and $g$ which is independent of the eigenvalue
    $\lambda$.   
  \item The functions $\phi_0, \phi_1, \ldots , \phi_{n_1}$ in
    (\ref{gensol}) that are arbitrary, depend only on the major variables
    $x^1_0,\; x^2_0, \ldots, x^m_0$, whereas the dependence on the minor
    variables of the functions $f_k$ and $g_k$ is restricted to certain
    fixed polynomials. For example, in the case when $m=1$, i.e., when the
    matrix consists of only one Jordan block, the functions $f_k$ and
    $g_k$ can be written as     
    \begin{align*}
      f_k(\xf) & = \lambda g_k(\xf) +
      k\sum_{j=0}^{k-1} P_{k-1,j}(x_1, x_2, \ldots , x_{k-j})
      \phi_k^{(j)}(x_0)\\  
      g_k(\xf) & = \sum_{j=0}^k P_{k,j}(x_1, x_2, \ldots ,
      x_{k-j+1})\phi_k^{(j)}(x_0),     
    \end{align*}
    where $\phi_k^{(j)}$ denote the $j$-th derivative of the function
    $\phi_k$, and $P_{k,j}$ are fixed polynomials, related to the
    partial Bell polynomials \cite{comtet}.    
\end{enumerate}

We provide a concrete example in order to illustrate how to read this
result from \cite{jodeit-1990}.  
\begin{Exmp3}\label{exn5}
  Consider the equation $($\ref{fmg}$)$ where $M$ is a constant
  $5\times 5$ matrix. As mentioned above, by changing independent
  variables in a linear way, we can assume that $M$ is in canonical
  Jordan form. We consider the equation $($\ref{fmg}$)$ in the
  particular case when $M$ is given by 
  \begin{align*}
    M = \left[
      \begin{array}{ccccc}
        \lambda_1 & 1 & 0 & 0 & 0\\
        0 & \lambda_1 & 0 & 0 & 0\\
        0 & 0 & \lambda_2 & 1 & 0\\
        0 & 0 & 0 & \lambda_2 & 0\\
        0 & 0 & 0 & 0 & \lambda_2\\
      \end{array}
    \right],
  \end{align*}
  where $\lambda_1$ and $\lambda_2$ are distinct eigenvalues. First of
  all, the general solution can be decomposed as $f = f^1 + f^2,$ $g =
  g^1 + g^2$, where $(f^1,\;g^1)$ and $(f^2,\;g^2)$ are solutions of
  \begin{align*}
    \nabla f^1 = \matris{\lambda_1}{1}{0}{\lambda_1}\nabla g^1,\quad and
    \quad \nabla f^2 =
    \matrist{\lambda_2}{1}{0}{0}{\lambda_2}{0}{0}{0}{\lambda_2} \nabla
    g^2,  
  \end{align*}
  respectively. From the formula $($\ref{gensol}$)$, we then obtain
  \begin{align*}
    \left\{ 
      \begin{aligned}
        f^1 & = f^1_0 + f^1_1 = \lambda_1g^1 + \phi_1(x_0) + c \\
        g^1 & = g^1_0 + g^1_1 = x_1\phi_1'(x_0) + \phi_0(x_0)
      \end{aligned}
    \right., 
  \end{align*}
  and
  \begin{align*}
    \left\{ 
      \begin{aligned}
        f^2 & = f^2_0 + f^2_1 = \lambda_2g^2 + \psi_1(y_0^1) + c \\
        g^2 & = g^2_0 + g^2_1 = \psi_0(y_0^1,\; y_0^2) + y^1_1\psi_1'(y_0^1)
      \end{aligned}
    \right.,
  \end{align*}
  where $x_0,\; y_0^1,\; y_0^2$ are the major variables, $x_1,\;
  y^1_1$ the minor variables, and $\phi_0,\; \phi_1,\; \psi_0,\; \psi_1$ are
  arbitrary scalar-valued analytic functions.
\end{Exmp3}

\subsection{The real case}\label{sec.gensolreal}

We consider now the equation (\ref{fmg}) in a convex domain of a
vector space over the real numbers. Since the field of real numbers is
not algebraically closed, the matrix $M$ is in general not similar to
a matrix of Jordan canonical form. Instead we assume that $M =
\mathrm{diag}(M^1,\ldots, M^p)$ has the real Jordan canonical form, so
that the blocks corresponding to the real eigenvalues have the same
structure as in the complex case, whereas a block $M^k =
\mathrm{diag}(M^k_1,\ldots, M^k_{r_k})$ corresponding to a complex
conjugated pair of eigenvalues, $\lambda, \bar{\lambda} = \alpha \pm
\mathrm{i}\beta$, consists of blocks $M^k_i$ of the form
\begin{align}\label{realjordan}
  \left[ 
    \begin{array}{cccc}
      \Lambda & I & &\\
      & \Lambda & \ddots &\\
      & & \ddots & I\\
      & & & \Lambda
    \end{array} \right],
  \quad {\rm where} \quad \Lambda =
  \matris{\alpha}{\beta}{-\beta}{\alpha}, \quad I =
  \matris{1}{0}{0}{1}.
\end{align} 
The general solution of (\ref{fmg}) can again be decomposed as in
(\ref{sumfg}), where each pair $(f^k,g^k)$ is a solution of the
subsystem $\nabla f = M^k \nabla g$. Thus, since the equations
corresponding to real eigenvalues can be treated in the same way as in
the complex case, it is enough to consider the case when $M =
\mathrm{diag}(M_1,\ldots, M_m)$ has only one pair of complex conjugate
eigenvalues $\lambda, \bar{\lambda}$. We assume that each $M_i$ is a
matrix of the form (\ref{realjordan}) of size $2(n_i +1)$ where $n_1
\ge \cdots \ge n_m \ge 0$. Moreover, we denote the variables
corresponding to the Jordan decomposition by $x^1_0, y^1_0, \ldots,
x^1_{n_1}, y^1_{n_1}, x^2_0, y^2_0, \ldots, x^m_{n_m}, y^m_{n_m}$, and
define complex variables $z^k_j = x^k_j + \mathrm{i}y^k_j$,
$\bar{z}^k_j = x^k_j - \mathrm{i}y^k_j$. In a similar way as in the
complex case, we introduce functions
\begin{align*}
  z^k_\mu = z^k_0 + z^k_1\mu + \cdots + z^k_{n_k}\mu^{n_k}, \quad k =
  1, 2, \ldots, m,
\end{align*}
depending on a real parameter $\mu$, and vectors
\begin{align*}
  \mathbf{z}^{(k)}_\mu = [z^1_\mu \quad z^2_\mu \;\; \cdots \;\;
  z^{\nu_k}_\mu]^T,  \quad k = 0, 1, \ldots, n_1,
\end{align*}
where $\nu_k$ is defined in the same way as above. If we let $F = f -
\bar{\lambda} g$, the general analytic solution of $\nabla f = M\nabla
g$ can be expressed as $F = F_0 + F_1 + \cdots + F_{n_1}$ where 
\begin{align}\label{Fkreal}
  F_k = \left. \left( \frac{\partial^k}{\partial \mu^k}\phi_k (
    \mathbf{z}^{(k)}_\mu ) - \sum_{l=1}^k \left(
      \frac{-\mathrm{i}}{2\beta} \right)^l \frac{k!}{(k-l)!}
    \overline{\frac{\partial^{k-l}}{\partial \mu^{k-l}}\phi_k (
      \mathbf{z}^{(k)}_\mu )}\right) \right|_{\mu=0}. 
\end{align}
and $\phi_k(s_1, \ldots, s_{\nu_k})$ are again arbitrary
complex-valued analytic functions depending on $\nu_k$ complex
variables. 

\subsection{A characterization of the general solution in terms of
  differentiable functions on algebras} 

In \cite{waterhouse-1992}, the general solution of (\ref{fmg}) is
characterized through components of functions which are differentiable
over algebras. From this characterization, some results regarding
multiplication of solutions can be obtained. However, the results in
\cite{waterhouse-1992}, regarding the equation (\ref{fmg}) is given
only for a quite restricted class of matrices $M$. We will give a
brief summary of the concepts and results in \cite{waterhouse-1992}
about the equation (\ref{fmg}).

Let $A$ be a commutative algebra of finite dimension over the real or
complex numbers, and consider a function $f:A\rightarrow A$. We say
that $f$ is \emph{$A$-differentiable} in the point $\af\in A$ if the
limit
\begin{align*}
\lim_{\substack{\xf \rightarrow 0 \\ \xf\in A^*}} \xf^{-1}(f(\af+\xf) - f(\af)) 
\end{align*}  
exists. When $A$ is a $\mathbb{C}$-algebra, any $A$-differentiable
function is analytic, i.e., it can be expressed locally as a power
series in the $A$-variable $\xf$. A function $V:A \rightarrow \R$ (or
$\mathbb{C}$) is called a \emph{component function} of an
$A$-differentiable function $f$ if it can be written as $V =
\phi \circ f$ for some linear function $\phi: A \rightarrow \R$ (or
$\mathbb{C}$).

When a $n\times n$ matrix $M$ consists of one single Jordan block,
there exists an algebra $A$ which is generated by a single element
$\mathbf{s}$ such that in some basis $\af_1, \ldots, \af_n$, $M$ is
the matrix of multiplication by $\mathbf{s}$, $\mathbf{s}\af_j =
\sum_i M_{ji} \af_i$. If $V = \phi (f)$ is a component function of a
$\mathcal{C}^2$ $A$-differentiable function $f$ and $W = \phi
(\mathbf{s}f)$ is the corresponding component function of
$\mathbf{s}f$, then $\nabla W = M\nabla V$. Conversely, if $V$ and $W$
are $\mathcal{C}^2$ functions such that $\nabla W = M\nabla V$, then
$V = \phi (f)$ and $W = \phi (\mathbf{s}f)$ for some
$A$-differentiable function $f$ and some linear function $\phi$.

This characterization in terms of component functions of
$A$-differentiable functions implies the existence of a multiplication
of solutions for the equation (\ref{fmg}) since products of
$A$-differentiable functions are again $A$-differentiable. Moreover,
in the complex case, since any $A$-differentiable function is
analytic, one can prove that any analytic solution of (\ref{fmg}) can
be expressed as a power series of a simple solution. Thus, since some
of the problems we are considering in this paper are already treated
in \cite{waterhouse-1992}, it is worth to emphasize in which way we
extend these results.

\begin{enumerate}
  \item In \cite{waterhouse-1992}, the only case considered is when
    the matrix $M$ has exactly one Jordan block corresponding to each
    eigenvalue. In this paper, the equation (\ref{fmg}) is studied for
    all matrices $M$.
  \item In \cite{waterhouse-1992}, only the existence of the algebra $A$, from the characterization of the general solution of (\ref{fmg}), is given. An explicit construction is natural and simple in the complex case, but not trivial in the real case. In this paper we provide an explicit multiplication formula for solutions of (\ref{fmg}) for every possible $M$, both in the real and in the complex case.
  \item The characterization of $A$-differentiable functions in terms of
    analytic functions is only valid for the complex case. We show that
    any solution of (\ref{fmg}), in both the complex and real case, can be
    expressed as a $*$-power series.  
  \item The equations (\ref{fmg}) constitute only a simple special case of a quite large family of systems of PDEs with $*$-multiplication. Therefore, the results obtained in this paper, for the restricted class of systems (\ref{fmg}), indicate that similar results may be established for the entire class of systems with $*$-multiplication.   
\end{enumerate}

In this paper, we have used the notation and conventions used in
\cite{jodeit-1990} and \cite{jonasson-2007}, rather than in
\cite{waterhouse-1992}. However, for the special cases when a
multiplication of solutions can be obtained from
\cite{waterhouse-1992}, this multiplication coincides, even though it 
may not be obvious immediately, with the $*$-multiplication studied
in this paper.

\section{Multiplication for systems $\nabla f = M\nabla
  g$}\label{cstmultsec} 

In this section we show that every equation (\ref{fmg}) admits a
$*$-multiplication. The approach we use for proving this, is to extend
the equation (\ref{fmg}) by introducing certain auxiliary dependent
variables and adding some equations that are consequences of the
original equation. One can then show that the extended system has the
form (\ref{avz}) and that it satisfies (\ref{azz}), so that it admits a
$*$-multiplication. Since there is a simple correspondence between
solutions of the original and the extended system, one can interpret
the $*$-multiplication as an algebraically defined operation on the
solution space of the equation (\ref{fmg}). We treat the real and
complex case separately, starting with the latter.

\begin{Rem3}
  We will sometimes use the exterior differential operator $\ed$
  instead of the gradient operator $\nabla$.  Since we are only
  considering complex analytic functions, the exterior differential
  reduces to the $\partial$-operator \cite{hormander}. Thus, we can
  write the equation $($\ref{fmg}$)$ as $\ed f = M^T \ed g$ which
  should be understood as the following equation for row matrices: 
  \begin{align*} 
    [\partial_1 f \; \cdots \; \partial_n f] = [\partial_1 g \; \cdots
    \; \partial_n g]M^T,   
  \end{align*} 
  where $\partial_i = \partial/\partial x_i$ denotes the partial
  derivative with respect to $x_i$, and $M^T$ denotes the transpose of
  the matrix $M$. 
\end{Rem3}

\subsection{The complex case}

When the general solution of the equation (\ref{fmg}) is obtained in
\cite{jodeit-1990}, Jodeit and Olver split their investigation into
three cases, and we will use the same strategy in this paper. The
cases are:
\begin{enumerate}
  \item The matrix $M$ consists of a single Jordan block, i.e.,
    \begin{align}\label{jblock}
      M = \left[ 
        \begin{array}{cccc}
          \lambda & 1 & &\\
          & \lambda & \ddots &\\
          & & \ddots & 1\\
          & & & \lambda
        \end{array} \right].
    \end{align} 
  \item $M=\mathrm{diag}(J_1,\ldots,\;J_m)$ has only one eigenvalue
    but consists of several Jordan blocks corresponding to that
    eigenvalue. 
  \item $M$ is a general complex matrix.
\end{enumerate}
Obviously, the first case is a special case of the second case, which
in turn is a special case of the third case. However, treating the
first two cases separately, has a didactic advantage.

\subsubsection{One Jordan block}

We assume now that the matrix $M$ consists of one Jordan block of size
$n+1$. By a change of dependent variable $h = f - \lambda g,$ the
equation (\ref{fmg}) reduces to the equation 
\begin{align}\label{hug}
  \nabla h = U_n \nabla g,
\end{align}
where $U_n$ is defined as the  the Jordan block of size $(n+1) \times
(n+1)$ corresponding to the eigenvalue $\lambda = 0$, i.e.,  
\begin{align}\label{un}
  U_n := \left[ 
    \begin{array}{cccc}
      0 & 1 & & \\
      & \ddots & \ddots & \\
      & & 0 & 1\\
      & & & 0
    \end{array} \right].
\end{align}
Equation (\ref{hug}) admits a $*$-multiplication on the solution
space: 
\begin{Prop3}\label{nil}
  Let $B$ be any nilpotent constant matrix. Then the system $\nabla
  f = B\nabla g$ can be extended to the system 
  \begin{align}\label{amur}
    A_\mu\nabla V_\mu \equiv 0 \quad \left(mod\; \mu^r \right), 
  \end{align}
  where 
  \begin{align*}
    A_\mu = -B + \mu I,\quad V_\mu = V_0 + \mu V_1 + \cdots +
    \mu^{r-3}V_{r-3} + \mu^{r-2}f + \mu^{r-1}g,  
  \end{align*}
  and $r$ is a natural number such that $B^r = 0$. The functions
  $V_0, \ldots ,V_{r-3}$ are auxiliary functions that are uniquely
  determined (up to additive constants) by any given solution $(f,g)$ 
  of $($\ref{hug}$)$.
\end{Prop3}
The proof relies on the following lemma.
\begin{Lem3}
  Let $B$ be a constant $n\times n$ matrix, and $f$ any smooth
  function such that $\ed (B\ed f) = 0$. Then $\ed (B^k\ed f) = 0$ for
  any $k\in\mathbb{N}$. 
\end{Lem3}
\begin{proof}
  The proof is by induction over $k$. For $k=2,$ by the assumption
  $\ed (B\ed f) = 0$, we have
  \begin{align*}
    \left( \ed (B^2 \ed f) \right)_{ij} & = \sum_{a,b=1}^n \left(
      B_{ab}B_{bi}\partial_a\partial_j f -
      B_{ab}B_{bj}\partial_a\partial_i f \right) \\   
    & = \sum_{a,b=1}^n \left(  B_{aj}B_{bi}\partial_a\partial_b f -
      B_{ai}B_{bj}\partial_a\partial_b f \right) = 0,  
  \end{align*}
  where $B_{ij}$ denotes the entry in row $i$ and column $j$ of the
  matrix $B$. Assume now that the statement is true for $k = p - 1$,
  where $p \ge 3$. If $\ed (B\ed f) = 0$, there exists, according to the
  Poincaré lemma, a function $g$ such that $\ed g = B \ed f$. Thus, we
  have $\ed (B \ed g) = \ed (B^2 \ed f) = 0$, and therefore by the
  inductive assumption
  \begin{align*}
    0 = \ed (B^{p-1} \ed g) = \ed (B^p \ed f).
  \end{align*} 
\end{proof}
\begin{proof}[Proof (of proposition \ref{nil}).]
  Let $A_\mu = -B + \mu I$, and $V_\mu = V_0 + \mu V_1 + \cdots +
  \mu^{r-1}V_{r-1}$. Then, the components of equation (\ref{amur})
  read 
  \begin{align*}
    \left\{ 
      \begin{aligned}
        0 & = B^T\ed V_0\\
        \ed V_0 & = B^Td V_1\\
        & \cdots\\
        \ed V_{r-2} & = B^Td V_{r-1},
      \end{aligned}
    \right.
  \end{align*}
  or equivalently (using the assumption that $B$ is nilpotent),
  \begin{align}\label{vbv}
    \ed V_{i} = \left( B^T \right)^{r-1-i}\ed V_{r-1},\quad i =
    0,1,\ldots ,r-2.
  \end{align}
  It is then obvious that if $V_\mu$ solves (\ref{vbv}), then the
  functions $f = V_{r-2}$, $g = V_{r-1}$ solve the equation $\nabla f
  = B\nabla g$. 

  On the other hand, assume now that $(f,g)$ solve the equation
  $\nabla f = B\nabla g$, and let $V_{r-2} = f$, $V_{r-1} = g$. Then,
  since $\ed (B^T \ed V_{r-1}) = \ed^2 V_{r-2} = 0$, the previous lemma
  implies that $\ed ((B^T)^k \ed V_{r-1}) = 0$, for any $k\ge 1.$
  Therefore, according to the Poincaré lemma there exist functions
  $V_0,V_1,\ldots ,V_{r-3}$ such that (\ref{vbv}) is satisfied.
\end{proof}

Since $A_\mu\nabla\mu^r \equiv 0$ is trivially satisfied, the equation
(\ref{amur}) admits a $*$-multi\-plica\-tion. Given any two solutions
$V_\mu$ and $\vt_\mu$ of (\ref{amur}), their $*$-product is given by
the formula
\begin{align*}
  V_\mu * \vt_\mu = \sum_{k=0}^{r-1}\left( \sum_{j=0}^k V_j\vt_{k-j}
  \right)\mu^k,
\end{align*}
or using the matrix notation 
\begin{align}\label{mult}
  \left[ \begin{array}{c} V_0 \\ V_1 \\ \vdots \\
      V_{r-1}\end{array}\right]*
  \left[ \begin{array}{c} \vt_0 \\ \vt_1 \\ \vdots \\
      \vt_{r-1}\end{array}\right] =
  \left[ \begin{array}{c} V_0\vt_0 \\ V_0\vt_1 + V_1\vt_0 \\ \vdots \\
      V_0\vt_{r-1} + V_1\vt_{r-2} + \cdots +V_{r-1}\vt_0
    \end{array}\right] 
\end{align}
Especially, since the matrix $U_n$ is nilpotent with $U_n^{n+1} = 0$,
proposition \ref{nil} allows us to extend the system (\ref{hug}) to
the $\mu$-dependent system 
\begin{align}\label{umu}
  (-U_n + \mu I)\nabla V_\mu \equiv 0 \quad \left(mod\; \mu^{n+1}
  \right).   
\end{align}
Thus, given two solutions $(h,g)$ and $(\tilde{h},\tilde{g})$ of
(\ref{hug}), there exist unique (up to additive constants) functions
$V_0,\ldots V_{n}$ and $\vt_0,\ldots \vt_{n}$ such that
\begin{align*}
  \nabla V_{i} = U_n^{n-i}\nabla g,\quad \textrm{and}\quad \nabla \vt_{i} =
  U_n^{n-i}\nabla \tilde{g},\quad i = 0,1,\ldots ,n.
\end{align*} 
A new solution of (\ref{hug}) is then algebraically constructed
through the $*$-multi\-plica\-tion as 
\begin{align*}
  (h,g) = \left( \sum_{i=0}^{n-1} V_{i}\vt_{n-1-i},\; \sum_{i=0}^n
    V_{i}\vt_{n-i} \right).
\end{align*}

\subsubsection{Several Jordan blocks corresponding to one eigenvalue}

Assume now that $M = \mathrm{diag}(J_1,\ldots ,J_m)$ has only one
eigenvalue $\lambda$, where each $J_i$ is a
Jordan block of size $n_i+1$, and $n_1\ge
n_2\ge\cdots \ge n_m$. Just as in the case of one single Jordan block,
by changing one dependent variable $h = f -\lambda g$, we obtain the
equivalent equation 
\begin{align}\label{hug2}
  \nabla h = U\nabla g, \quad \textrm{where}\quad U =
  \mathrm{diag}(U_{n_1},\ldots, U_{n_m}), 
\end{align}
and $U_i$ is given by (\ref{un}). Since $U$ is nilpotent, we can apply
proposition \ref{nil}, with $B = U$ and $r = n_1 + 1$, in order to
extend (\ref{hug2}) to the system
\begin{align}\label{umu2}
  (-U+\mu I)\nabla V_\mu \equiv 0 \quad \left(mod\; \mu^{n_1+1}
  \right),  
\end{align}
Thus, there is a $*$-multiplication for solutions of the equation
(\ref{fmg}) (having the same multiplication formula (\ref{mult})) also 
for the case when $M$ consists of several Jordan blocks corresponding
to the same eigenvalue.
 
Besides generating solutions of the system (\ref{hug2}) with
$*$-multiplication, it follows from the next proposition that one can
also embed solutions of the subsystems $\nabla h =
\mathrm{diag}(U_{n_i}, \ldots, U_{n_m})\nabla g$, for any $1\le i \le
m$, into the solution space of (\ref{hug2}).

\begin{Prop3}\label{embedd} Suppose that $V_\mu$ is a solution of the
  equation
  \begin{align}\label{amur2} A_\mu\nabla V_\mu \equiv 0 \quad
    \left(mod\; \mu^r \right),\quad V_\mu = \sum_{i=0}^{r-1}V_i\mu^i
  \end{align}
  defined in a subset $\Omega_1$ of a vector space $\mathcal{V}_1$,
  where $A_\mu = -A + \mu I$, and $r$ is a natural number. Then $W_\mu
  = \mu^s V_\mu$ is a solution of the extended system   
  \begin{align*}
    B_\mu\nabla W_\mu \equiv 0 \quad \left(mod\; \mu^{r+s} \right),\quad
    W_\mu = \sum_{i=0}^{r+s-1}W_i\mu^i         
  \end{align*}
  defined in $\Omega_1\times \mathcal{V}_2,$ where $\mathcal{V}_2$ is 
  another vector space, $B_\mu = B + \mu I$, and $B$ is the square
  matrix with the block structure
  \begin{align*}
    B = -\matris{C}{0}{0}{A},     
  \end{align*}
  for some constant matrix $C$ and some $s\in \mathbb{N}$.
\end{Prop3}

\begin{proof}
  Let $V_\mu$ be a solution of the equation (\ref{amur2}), i.e., there
  exists a $1$-form $\alpha$, not depending on $\mu$, such that
  $A_\mu^T\ed V_\mu = \mu^r\alpha$. Then
  \begin{align*}
    B_\mu^T\ed (\mu^sV_\mu) & = \mu^s B_\mu^T\ed V_\mu = \mu^s \left(
      C_\mu^T\ed_1 V_\mu + A_\mu^T\ed_2 V_\mu\right)\\
    & = \mu^s A_\mu^T\ed_2 V_\mu = \mu^{s+r}\alpha \equiv 0 \quad
    \left(mod\; \mu^{r+s} \right), 
  \end{align*}
  where $\ed_1$ and $\ed_2$ denote the exterior differential operators
  on $\mathcal{V}_1$ and $\mathcal{V}_2$, respectively.
\end{proof}

Let $(h^i,g^i)$ be a solution of $\nabla h^i =\mathrm{diag}(U_{n_i},\ldots,
U_{n_m})\nabla g^i$. Then, according to proposition \ref{nil}, there exists
a unique (up to unessential constants) solution $V^i_\mu = V_0^i + \cdots
+ V^i_{n_i-2}\mu^{n_i-2} + h^i\mu^{n_i-1} + g^i\mu^{n_i}$ of the extended
system 
\begin{align*}
  \left( -\mathrm{diag}(U_i,\ldots,U_m) + \mu I \right) \ed V^i_\mu
  \equiv 0 \quad \left( mod\; \mu^{n_i+1} \right).     
\end{align*}
According to proposition \ref{embedd}, it then follows that $V_\mu^i$
can be embedded as a solution $\mu^{n_1-n_i}V_\mu^i$ of the equation
(\ref{hug2}).   

\subsubsection{The general complex case}

Since the general solution of equation (\ref{fmg}) can be written
as a sum (\ref{sumfg}) where $(f^i,g^i)$ is obtained from the general
solution of the system (\ref{hug2}), it is clear that there is no
non-trivial $*$-multiplication for the system (\ref{fmg}) other than
the multiplications induced from the subsystems corresponding to
single eigenvalues, that was introduced above.

In order to illustrate the $*$-multiplication for equations of the
form (\ref{fmg}), we return to example \ref{exn5}.

\begin{Exmp3}
  Consider again equation $($\ref{fmg}$)$ in the case when $M$ is the
  $5\times 5$ matrix defined in example \ref{exn5}. Take for instance the
  particular solution
  \begin{align*}
    \left\{
      \begin{aligned}
        f & = \lambda_1 2x_0x_1 + \lambda_1(x_0)^3 + (x_0)^2 +
        \lambda_2 y_0^1(y_0^2)^2 + \lambda_2 y_1^1 + y_0^1\\
        g & = 2x_0x_1 + (x_0)^3 + y_0^1(y_0^2)^2 + y_1^1
      \end{aligned}
    \right.
  \end{align*}
  This solution can be written as a sum of $*$-products of simple
  (first order polynomials) solutions of subsystems. In detail, 
  \begin{align*}
    \left\{
      \begin{aligned}
        f & = \lambda_1g_1 + h_1 + \lambda_2g_2 + h_2 \\
        g & = g_1 + g_2
      \end{aligned}
    \right.
  \end{align*}
  where
  \begin{align*}
    h_1 + \mu g_1 = \mu * (x_0 + \mu x_1)^3_* + (x_0 + \mu x_1)^2_*,
  \end{align*}
  and
  \begin{align*}
    h_2 + \mu g_2 = \left( \mu (y_0^2)^2_* \right)*(y_0^1 + \mu y_1^1)
    + (y_0^1 + \mu y_1^1). 
  \end{align*}
  We note that the factor $\mu (y_0^2)^2_*$ is the embedding,
  according to proposition \ref{embedd}, of the solution $(y_0^2)^2_*
  = (y_0^2)^2$ of the subsystem $\nabla h = U_0 \nabla g$ (or
  equivalently $h'(y_0^2) = 0$), into the system $\nabla h_2 =
  \mathrm{diag}(U_1,U_0)\nabla g_2$.  
\end{Exmp3}

The possibility to express solutions of an equation of the kind
(\ref{fmg}) as a sum of simple $*$-products, that was illustrated in
the previous example, is a general property. In the next section we
show that any analytic solution of (\ref{fmg}) can be expressed as a
power series of simple solutions with respect to the
$*$-multiplication.  

\subsection{The real case}

We consider now the equation (\ref{fmg}) over a real vector space when
the matrix $M$ has real entries. From the discussion in section
\ref{sec.gensolreal}, it is clear that it is no restriction to assume
that $M$ has only one pair of complex conjugate eigenvalues. Thus, we
let $M = \mathrm{diag}(M_1,\ldots, M_m)$ be an arbitrary square matrix
having real Jordan form with eigenvalues $\lambda, \bar{\lambda} =
\alpha \pm \mathrm{i}\beta$, and we use the same notation as in
section \ref{sec.gensolreal}. Furthermore, we can, without loss of
generality, assume that $\alpha = 0$ and $\beta = 1$. This is realized
by changing both dependent and independent variables according to
\begin{align*}
  \tilde{f} & = f - \alpha g, \quad \tilde{g} = \beta g, \\
  \tilde{x}_j^k & = \beta^{j-n_k}x_j^k, \quad \tilde{y}_j^k =
  \beta^{j-n_k}y_j^k, \quad j = 0, 1, \ldots, n_k \quad k = 1, 2,
  \ldots, m.    
\end{align*}

The matrix $M$ is not nilpotent, which was the case for the complex
case, but the system (\ref{fmg}) is a  so called quasi-Cauchy--Riemann
equation which are known to admit $*$-multiplication. 

\begin{Thm3}\label{realmult}
  Let $M$ be a square matrix of size $2(n + 1)$ having real Jordan
  form, and with no other eigenvalues than $\pm \mathrm{i}$. Then
  equation $($\ref{fmg}$)$ can be extended to the system    
  \begin{align}\label{amuv}
    A_\mu\nabla V_\mu \equiv 0 \quad \left(mod\; (\mu^2 + 1)^{n+1} \right), 
  \end{align}
  where 
  \begin{align*}
    A_\mu = M^{-1} + \mu I,\quad V_\mu = V_0 + V_1 \mu + \cdots +
    V_{2n+1}\mu^{2n+1},  
  \end{align*}
  and $V_0 = f$ and $V_{2n+1} = g$. The functions $V_1, \ldots, V_{2n}$
  are uniquely determined (up to additive constants) by the solution
  $(f,g)$.  
\end{Thm3} 

\begin{proof}
  Since $\det{M} = (\det{\Lambda})^{n+1} = 1$, equation (\ref{fmg}) can
  be written as 
  \begin{align}\label{mfdetmg}
    M^{-1}\nabla f = \det{M^{-1}}\nabla g.
  \end{align}
  Equation (\ref{mfdetmg}) is a quasi-Cauchy--Riemann equation and is
  therefore equivalent to the parameter dependent equation (see
  \cite{jonasson-2006})\footnote{In \cite{jonasson-2006}, the
    definition of a quasi-Cauchy--Riemann equation require that the
    matrix $M$ is symmetric, but as long as the entries are constant
    this assumption is not necessary.}:
  \begin{align}\label{fdetmu}
    (M^{-1} + \mu I)\nabla V_\mu \equiv 0 \quad \left(mod\;
      \det{(M^{-1} + \mu I)} \right), 
  \end{align} 
  where $V_0 = f$ and $V_{2n+1} = g$. Thus, since 
  \begin{align*}
    \det{(M^{-1} + \mu I)} = (\det{(\Lambda + \mu I)})^{n+1} = (\mu^2
    + 1)^{n+1}, 
  \end{align*}
  the proof is complete. 
\end{proof}

Since $(M^{-1} + \mu I)\nabla (\mu^2 + 1)^{n+1} = 0$, the solution
space of the system (\ref{amuv}) admits $*$-multiplication. In the
special case when (\ref{fmg}) is the Cauchy--Riemann equations, i.e.,
$M = \Lambda$ with $\alpha = 0$ and $\beta = 1$, the parameter
dependent equation (\ref{fdetmu}) coincides with the original equation
(\ref{fmg}) and the $*$-multiplication reduces to the ordinary
multiplication, obtained from the multiplication of holomorphic
functions.

In a similar way as in the complex case, solutions of (\ref{fmg}) can
be constructed by embedding of solutions of subsystems. Let $V_\mu$ be
a solution of  
  \begin{align*}
    \tilde{A}_\mu\nabla V_\mu \equiv 0 \quad \left(mod\; (\mu^2 +
      1)^{\tilde{n}+1} \right),  
  \end{align*}
where $\tilde{A}_\mu = \mathrm{diag}(M_i^{-1},\ldots, M_m^{-1}) + \mu
I$ and $\tilde{n} = n_i + \cdots + n_m + m - i$. Then $(\mu^2 +
1)^{n-\tilde{n}}V_\mu$ is a solution of (\ref{amuv}). A proof of a
more general result is given in \cite{jonasson-2007}.  

\section{Analytic solutions are $*$-analytic}\label{pssec}

We have seen in the previous section that any equation of the form
(\ref{fmg}) can be reduced to a number of systems of the form
(\ref{umu2}) (or (\ref{amuv}) in the real case), each one equipped
with a $*$-multiplication on the solution space. In this section we
prove that every analytic solution can be expressed as a power series
of simple solutions, similarly as for the Cauchy--Riemann equations.
Being precise, by simple solution we mean a solution which is linear
in the independent variables. Furthermore, we say that a solution
$(f,g)$ of (\ref{fmg}) is $*$-\emph{analytic} if it can be expressed
locally as a finite sum of power series of simple solutions with
respect to the different $*$-multiplications of the corresponding
subsystems of the form (\ref{umu2}) and (\ref{amuv}). In other words,
in this section we aim to prove that every analytic solution of
(\ref{fmg}) is also $*$-analytic. Since it is no restriction, we will
only consider power series expansions in a neighborhood of the origin.

In both the real and the complex case, equation (\ref{fmg}) is
extended to a finite number of $\mu$-dependent equations of the form 
\begin{align}\label{xmuvmu}
  (X + \mu I) \nabla V_\mu \equiv 0 \quad \left(mod\; Z_\mu \right), 
\end{align} 
where $X$ is a constant matrix, $Z_\mu = Z_0 + \cdots + Z_n\mu^n$ is a 
polynomial with constant coefficients, $V_\mu = V_0 + \cdots +
V_n\mu^n$, and $V_n = g$. Equation (\ref{xmuvmu}) can be written
without the parameter $\mu$ as a system
\begin{align}\label{xvj}
  X \nabla V_j + \nabla V_{j-1} = Z_j \nabla g \quad j = 0, 1, \ldots,
  n,  
\end{align}
where $V_{-1} = 0$. Let $V_\mu^{(N)} = V_0^{(N)} + \cdots +
V_n^{(N)}\mu^n$ be a solution with coefficients being polynomials of 
degree $N$ in the coordinate functions,
\begin{align*}
  V_j^{(N)} = \sum_{k=0}^N \sum_{\substack{I=(i_0,\ldots,i_n) \in
      \mathbb{N}^n \\ i_0 + \cdots + i_n = k}} a_I x_0^{i_0}x_1^{i_1}
  \cdots x_n^{i_n}. 
\end{align*}
Because of the relation (\ref{xvj}), the power series 
\begin{align*}
  \lim_{N \rightarrow \infty} V_j^{(N)}, \quad  j = 0, 1, \ldots, n
\end{align*}
will all have the same domain of convergence. This means that when
constructing $*$-power series of a simple solution, it is enough to
consider the domain of convergence for the power series defined by the
coefficient at the highest power $\mu^n$. 

\subsection{The complex case}

We split the investigation of the complex case into three different
cases in the same way as in the previous section.

\subsubsection{One Jordan block}

We consider first the equation (\ref{fmg}) when the Jordan
decomposition of the constant matrix $M$ consist of one single Jordan
block. We have seen that, then there is no loss of generality to
consider the equation (\ref{hug}) instead.
\begin{Thm3}\label{gensolthm1}
  Let $(h,g)$ be an analytic solution of the equation $($\ref{hug}$)$.
  Then there exist constant solutions $(b_i)_\mu = b_{i0} + b_{i1}\mu
  + \cdots + b_{in}\mu^n$ (with $b_{ij} \in \mathbb{C}$) of the
  extended system $($\ref{umu}$)$ for $i = 0, 1, \ldots$, such that  
  \begin{align}\label{ps1}
    \sum_{i=0}^\infty (b_i)_\mu*(x_\mu)_*^i = V_0 + \cdots +
    V_{n-2}\mu^{n-2} + h\mu^{n-1} + g\mu^n,
  \end{align}
  where $x_\mu = x_0 + x_1\mu + \cdots +  x_n\mu^n$.
\end{Thm3}
We formulate a result about the structure of the $*$-powers
$(x_\mu)^i_*$ as a lemma: 
\begin{Lem3}\label{xi}
  For each $i\in \mathbb{N}$ we have
  \begin{align}\label{xif}
    (x_\mu)^i_* = \sum_{j=0}^n \Psi_{i,j}(\xf)\mu^j, 
  \end{align}
  where $\Psi_{i,j}$ are the polynomials given by
  \begin{align*}
    \Psi_{i,j}(\xf) & := \! \sum_{\substack{a_0,a_1,\ldots,a_n \in \mathbb{N}
        \\ a_0 + a_1 + \cdots + a_n = i \\ a_1 + 2a_2 + \cdots + na_n
        = j }}\! \binom{i}{a_0, a_1, \ldots, a_n}x_0^{a_0}x_1^{a_1}
    \cdots x_n^{a_n} \\ 
    & = \sum_{s=0}^j \binom{i}{s}x_0^{i\!-\!s} B_{js}(x_1, 2x_2, \ldots,
    (j\!-\!s\!+\!1)!x_{j\!-\!s\!+\!1}),
  \end{align*}
  where $B_{js} =B_{js}(x_1, x_2, \ldots, x_{j\!-\!s\!+\!1}) $ are the
  partial Bell polynomials \cite{comtet}  given by
  \begin{align*}
    \sum_{\substack{a_1,a_2,\ldots,a_{j\!-s\!+1} \in \mathbb{N}
        \\ a_0 + a_1 + \cdots + a_{j\!-s\!+1} = s \\ a_1 + 2a_2 +
        \cdots + (j\!-s\!+1)a_{j\!-s\!+1}
        = j }}\!\!\!\!\!\!\!\!\!\!\!\! \binom{s}{a_1, \ldots,
      a_{j\!-\!s\!+\!1}}\!\!\left(\frac{x_1}{1!}\right)^{a_1}
    \left(\frac{x_2}{2!}\right)^{a_2} \cdots  
    \left( \frac{x_{j\!-\!s\!+\!1}}{(j\!-\!s\!+\!1)!}
    \right)^{a_{j\!-\!s\!+\!1}}  
  \end{align*}
\end{Lem3}
\begin{proof}
  From the multinomial theorem, we have
  \begin{align*}
    x_\mu^i = \sum_{\substack{a_0,a_1,\ldots,a_n \in \mathbb{N}
        \\ a_0 + a_1 + \cdots + a_n = i }}\! \binom{i}{a_0, a_1,
      \ldots, a_n}x_0^{a_0}x_1^{a_1} \cdots x_n^{a_n}\mu^{a_1 \!+\!
      2a_2 \!+ \cdots +\! na_n}.
  \end{align*}
  Thus, formula (\ref{xif}) follows now from the definition of the
  $*$-multiplication. Moreover, 
  \begin{align*}
    \Psi_{i,j}(\xf) & = \sum_{a_0=0}^i\frac{x_0^{a_0}}{a_0!} \!
    \sum_{\substack{a_1,a_2,\ldots,a_n \in \mathbb{N} \\ a_1 + \cdots
        + a_n = i - a_0 \\ a_1 + 2a_2 + \cdots + na_n
        = j }}\! \frac{i!}{a_1!a_2!\cdots a_n!}x_1^{a_1} \cdots
    x_n^{a_n} \\  
    & = \sum_{s=0}^j \binom{i}{s}x_0^{i-s}\!\!\!\!\!\!\!\!\!
    \sum_{\substack{a_1,a_2,\ldots,a_{j\!-\!s\!+\!1} \in \mathbb{N} \\ a_1 + \cdots
        + a_{j\!-\!s\!+\!1} = s \\ a_1 + 2a_2 + \cdots + (j\!-\!s\!+\!1)a_{j\!-\!s\!+\!1}
        = j }}\! \binom{s}{a_1, \ldots, a_{j\!-\!s\!+\!1}}
    x_1^{a_1} \cdots x_{j\!-\!s\!+\!1}^{a_{j\!-\!s\!+\!1}}\\
    & = \sum_{s=0}^j \binom{i}{s}x_0^{i\!-\!s} B_{js}(x_1, 2x_2, \ldots,
    (j\!-\!s\!+\!1)!x_{j\!-\!s\!+\!1}). 
  \end{align*} 
\end{proof}

\begin{proof}[Proof (of theorem \ref{gensolthm1}).]
  For any solution $(h,g)$ of (\ref{hug}), $h$ is uniquely determined
  from $g$ (up to an additive constant). Thus, it is enough to prove
  that for any analytic solution $(h,g)$, we can choose $(b_i)_\mu$ such
  that the highest order coefficient of $\sum_{i=0}^\infty
  (b_i)_\mu*(x_\mu)^i_*$ coincides with $g$. From the results obtained
  in \cite{jodeit-1990}, presented in section \ref{fmgsec}, the general
  solution $g$ is given by
  \begin{align}\label{geng}
    g(\xf) = \left. \sum_{j=0}^n \frac{\partial^j}{\partial
        \mu^j}\phi_j \left( x_\mu \right) \right|_{\mu=0}, 
  \end{align}
  where $\phi_1, \phi_2, \ldots, \phi_n$ are arbitrary analytic
  functions of one complex variable, given by 
  \begin{align}\label{psphi}
    \phi_j(s) = \sum_{i=0}^\infty c_{ij}s^i, \quad j = 0, 1, \ldots,
    n.  
  \end{align}
  Hence, if we write the function $g$ in a more explicit way, we
  obtain
  \begin{align*}
    g(\xf) & = \sum_{j=0}^n\sum_{i=0}^\infty c_{ij} \left. \left(
        \frac{\partial^j}{\partial \mu^j} (x_\mu)^i \right)
    \right|_{\mu=0} \\
    & = \sum_{j=0}^n\sum_{i=0}^\infty c_{ij}\!\!\! \left. \left( \!\!
        \frac{\partial^j}{\partial \mu^j} \!\!\!\!
        \sum_{\substack{a_0,\ldots,a_n \in \mathbb{N} \\ a_0 + \cdots
            + a_n = i }}\!\!\!\! \binom{i}{a_0, \ldots,
          a_n}x_0^{a_0}x_1^{a_1} \cdots x_n^{a_n}\mu^{a_1 \!+\!  
          \cdots \!+\! na_n} \right) \!\right|_{\mu=0}\\
    & = \sum_{j=0}^n\sum_{i=0}^\infty j!c_{ij} \Psi_{i,j}.
  \end{align*}
  We define now the constant solutions $(b_i)_\mu$ as 
  \begin{align*}
    (b_i)_\mu = \sum_{k=0}^n (n-k)!c_{i,n-k}\mu^k.
  \end{align*}
  From lemma (\ref{xi}), it then follows that the highest order
  coefficient of the $*$-power series $V_\mu = \sum (b_i)_\mu *
  (x_\mu)_*^i$ coincides  with $g(\xf)$.  

\end{proof}

We illustrate this result with a simple example.

\begin{Exmp3}
  Consider the system $\nabla h = U_2\nabla g$, where $U_2$ is the
  $3\times 3$ matrix described above, i.e., 
  \begin{align}\label{u3ex}
    \nabla h = \matrist{0}{1}{0}{0}{0}{1}{0}{0}{0} \nabla g.
  \end{align}
  The general analytic solution of $($\ref{u3ex}$)$ is given by 
  \begin{align}\label{u3exsol}
    \left\{
      \begin{aligned}
        h & = \phi_1(x_0) + 2x_1\phi_2'(x_0) + c\\
        g & = \phi_0(x_0) + x_1\phi_1'(x_0) + 2x_2\phi_2'(x_0) +
        x_1^2\phi_2''(x_0) 
      \end{aligned}
    \right.
  \end{align} 
  where $\phi_0,\phi_1,\phi_2$ are general analytic functions of one
  complex variable. If the corresponding power series representations
  are given by $($\ref{psphi}$)$ for $j = 0,1,2$, then the general
  solution $($\ref{u3exsol}$)$ is $*$-analytic with the following power
  series representation:  
  \begin{align*}
    \vektort{V_0}{h}{g} = \sum_{i=0}^\infty
    \vektort{c_{i2}}{c_{i1}}{c_{i0}}*\vektort{x_0}{x_1}{x_2}_*^i.
  \end{align*}
\end{Exmp3}

\subsubsection{Several Jordan blocks corresponding to one eigenvalue}

Assume now that the matrix $M$ in equation (\ref{fmg}) consists of
only one eigenvalue $\lambda$, but that its Jordan canonical form
consists of several Jordan blocks. As we have seen, it is then no
restriction to assume that the eigenvalue is zero, so that (\ref{fmg})
is reduced to the equation (\ref{hug2}). Then a given analytic
solution can in general not be written directly as a power series of
simple solutions with respect to the $*$-multiplication (\ref{mult})
for that system. This is easily realized from the following example. 

\begin{Exmp3}
  Consider the equation $\nabla h = \mathrm{diag}(U_1,U_0) \nabla g$,
  i.e.,  
  \begin{align}\label{u1u0}
    \nabla h = \matrist{0}{1}{0}{0}{0}{0}{0}{0}{0} \nabla g.
  \end{align}
  The general analytic solution is given by 
  \begin{align*}
    \left\{
      \begin{aligned}
        h & = \phi(x) + c \\
        g & = y\phi'(x) + \psi(x,z) 
      \end{aligned}
    \right.
  \end{align*} 
  where $\phi$ and $\psi$ are arbitrary analytic functions, and the
  $*$-multiplication formula is given by
  \begin{align}\label{u1u0mult}
    \vektor{h}{g}*\vektor{\tilde{h}}{\tilde{g}} =
    \vektor{h\tilde{h}}{h\tilde{g} + g\tilde{h}}.
  \end{align}
  Thus, we see that the $*$-product of polynomial solutions will never
  have higher degree in the $z$-variable than the highest degree of
  $g$ and $\tilde{g}$. Since a polynomial solution can have any given
  polynomial degree in $z$, it becomes obvious that not every analytic
  solution of $($\ref{u1u0}$)$ can be written as a power series of
  simple solutions with respect to the $*$-multiplication
  $($\ref{u1u0mult}$)$.

  The remedy is to embed power series solutions of the trivial
  subsystem $\nabla h = U_0 \nabla g$ in the sense of proposition
  \ref{embedd}. The general solution is given by $g = \tau(z)$ where
  $\tau$ is an arbitrary function, and the $*$-multiplication
  coincides with the ordinary multiplication. Thus, according to
  proposition \ref{embedd}, we can embed $*$-powers $z_*^j = z^j$ of
  the trivial solution $g = z$ into solutions of the system
  $($\ref{u1u0}$)$. Any analytic solution of $($\ref{u1u0}$)$ can then
  be written as an infinite sum of products of the solutions $(x,
  y)^i_*$ and the embedded solutions $(0, z^j)$ of the $*$-powers
  $z_*^j$: 
  \begin{align*}
    \vektor{\displaystyle\sum_{i=0}^\infty
      c_ix^i}{\displaystyle\sum_{i=0}^\infty ic_iyx^i +
      \sum_{i,j=0}^\infty c_{ij}x^iz^j} = \sum_{i=0}^\infty
    c_i\vektor{x}{y}^i_* + \sum_{i=0}^\infty
    c_{ij}\vektor{x}{y}^i_**\vektor{0}{z_*^j}  
  \end{align*}
\end{Exmp3}

In the general case with one eigenvalue, we proceed in the same way as
in the example, and embed power series solutions for the
subsystems $\nabla h^i =\mathrm{diag}(U_{n_i},\ldots, U_{n_m})\nabla
g^i$, in the sense of proposition \ref{embedd}.
\begin{Thm3}\label{ps2}
  Every analytic solution of the matrix equation $($\ref{hug2}$)$ is
  $*$-analytic. In detail, let $(h,g)$ be an analytic solution defined
  by 
  \begin{align*}
    g = \sum_{j=0}^{n_1} \left. \frac{\partial^j}{\partial
        \mu^j}\phi_j \left( \xf^{(j)}_\mu\right) \right|_{\mu=0},  
  \end{align*}  
  where $\phi_0,\phi_1,\ldots ,\phi_{n_1}$ are analytic functions with
  power series representations
  \begin{align*}
    \phi_j(s_1,s_2,\ldots ,s_{\nu_j}) = \sum_{I=(i_1,i_2, \ldots ,
      i_{\nu_j})\in\mathbb{N}^{\nu_j}} c_{I,j}s_1^{i_1}s_2^{i_2} \cdots
    s_{\nu_j}^{i_{\nu_j}}. 
  \end{align*}
  Then $(h,g)$ is $*$-analytic with power series representation
  \begin{align*}
    V_0 + \cdots + V_{n_1-2}\mu^{n_1-2} + h\mu^{n_1-1} +
    g\mu^{n_1} = \sum_{r=1}^m \sum_{I=(i_1,i_2, \ldots ,
      i_r)\in\mathbb{N}^r} \!\!\!\!\!\!(b_I)_\mu * (\xf_\mu)^I_*,
  \end{align*}
  where, for each multi-index $I = (i_1,i_2, \ldots ,i_r)$,
  $(b_I)_\mu$ is the trivial solution of the extended system
  $($\ref{umu2}$)$ defined as 
  \begin{align*}
    (b_I)_\mu = \sum_{j=n_{r+1}+1}^{n_r} j!c_{I,j}\mu^{n_r-j}, \quad
    (n_{m+1} :=- 1),
  \end{align*}
  and
  \begin{align}\label{xI}
    (\xf_\mu)^I_* := \left( x^1_\mu \right)_*^{i_1} * \left(
      \mu^{n_1-n_2}(x^2_\mu)_*^{i_2} * \cdots * \left(
      \mu^{n_{r-1}-n_r}(x^r_\mu)_*^{i_r} \right)\right). 
  \end{align}
\end{Thm3}

\begin{Rem3}
  The definition $($\ref{xI}$)$ of the powers $(\xf_\mu)^I_*$ needs
  some further explanation since the $*$-operators on the right hand
  side of $($\ref{xI}$)$ denote in fact multiplications for different
  systems. In order to calculate $(\xf_\mu)^I_*$, one has to start with
  the power $(x^r_\mu)_*^{i_r}$, which is a solution of the subsystem
  $\nabla h = \mathrm{diag}(U_r, \ldots, U_m)\nabla g$. This solution is
  then embedded, in the sense of proposition \ref{embedd}, into the
  solution $\mu^{n_{r-1}-n_r}(x^r_\mu)_*^{i_r}$ of the subsystem $\nabla
  h = \mathrm{diag}(U_{r-1}, U_r, \ldots, U_m)\nabla g$, and is
  thereafter $*$-multiplied with $(x^{r-1}_\mu)_*^{i_{r-1}}$. One
  continues this successive embedding of solutions into larger
  subsystems until the solution $(x^2_\mu)_*^{i_2} * \cdots *
  (\mu^{n_{r-1}-n_r}(x^r_\mu)_*^{i_r})$ of the system $\nabla h =
  \mathrm{diag} (U_2, \ldots, U_m)\nabla g$ is embedded (through
  multiplication with $\mu^{n_1-n_2}$) into a solution of the system
  $($\ref{hug2}$)$, and then finally multiplied with the solution
  $(x^1_\mu)_*^{i_1}$.
\end{Rem3}

As in the case with only one Jordan block, we place the calculation of
the powers $(\xf_\mu)^I_*$ in a separate lemma. 

\begin{Lem3}
  For each $I = (i_1,i_2, \ldots ,i_r)\in \mathbb{N}^r$ we have
  \begin{align*}
    (\xf_\mu)^I_* =
    \sum_{j=0}^{n_r} \Psi_{I,j}(\xf^1, \xf^2, \ldots,
    \xf^r)\mu^{n_1-n_r+j},  
  \end{align*}
  where $\Psi_{I,j}$ are polynomials defined through the polynomials
  $\Psi_{i,j}$ in lemma \ref{xi} as  
  \begin{align*}
    \Psi_{I,j}(\xf^1, \xf^2, \ldots, \xf^r) = \sum_{\substack{ 0\le
        j_1,\ldots, j_r \le n_r\\
        j_1+j_2+\cdots
      +j_r=j}} \Psi_{i_1,j_1}(\xf^1)\Psi_{i_2,j_2}(\xf^2) \cdots
    \Psi_{i_r,j_r}(\xf^r)  
  \end{align*}
\end{Lem3}

\begin{proof}
  For $r = 1$, the statement reduces to lemma \ref{xi}. For arbitrary
  $r$, since 
  \begin{align*}
    (\xf_\mu)^I_* = (x^1_\mu)^{i_1}_* * \left(
      (\tilde{\xf}_\mu)^{(i_2, \ldots i_r)}_* \mu^{n_1-n_2}\right),   
  \end{align*}
  where $\tilde{\xf} = [\xf^2, \ldots, \xf^r]^T$, we may assume by
  induction that 
  \begin{align*}
    (\xf_\mu)^I_* & = \left( \sum_{s=0}^{n_1} \Psi_{i_1,s}(\xf^1)\mu^s
    \right) \! * \!\left( \mu^{n_1\!-n_2}\!\sum_{t=0}^{n_r} \Psi_{(i_2, \ldots
        i_r),t}(\xf^2, \ldots, \xf^r)\mu^{n_2-n_r+t} \right). 
  \end{align*} 
  Thus, from the $*$-multiplication formula (\ref{mult}), we obtain 
  \begin{align*}
    (\xf_\mu)^I_* & = \sum_{j=0}^{n_r} \left( \sum_{s=0}^j
      \Psi_{i_1,s}(\xf^1)\Psi_{(i_2, \ldots, i_r),j-s}(\xf^2, \ldots,
      \xf^r)\right)\mu^{n_1-n_r+j}\\
    & = \sum_{j=0}^{n_r}\Psi_{I,j}(\xf^1, \xf^2, \ldots,
    \xf^r)\mu^{n_1-n_r+j},   
  \end{align*}
  where the last equality follows immediately from the definition of
  the polynomials $\Psi_{I,j}$.
\end{proof}

\begin{proof}[Proof (of theorem \ref{ps2}).]
  Since $\nu_{n_{r+1}+1} = \nu_{n_{r+1}+2} = \cdots = \nu_{n_r} = r$
  for $r = 1, 2, \ldots, m$, we have 
  \begin{align*}
    g & = \sum_{j=0}^{n_1} \sum_I \left. c_{I,j}
      \frac{\partial^j}{\partial \mu^j} \left(
        (x^1_\mu)^{i_1}(x^2_\mu)^{i_2} \cdots
        (x_\mu^{\nu_j})^{i_{\nu_j}} \right) \right|_{\mu=0}\\
    & = \sum_{r=1}^{m} \sum_{j=n_{r+1}+1}^{n_r} \sum_I \left. c_{I,j}
    \frac{\partial^j}{\partial \mu^j} \left(
      (x^1_\mu)^{i_1}(x^2_\mu)^{i_2} \cdots (x_\mu^r)^{i_r} \right)
  \right|_{\mu=0}. 
  \end{align*}
  Thus, it is enough to show that for each $I = (i_1,i_2, \ldots
  ,i_r)$, the expression 
  \begin{align*}
    \sum_{j=n_{r+1}+1}^{n_r} \left. c_{I,j} \frac{\partial^j}{\partial
        \mu^j} \left( (x^1_\mu)^{i_1}(x^2_\mu)^{i_2} \cdots
        (x_\mu^r)^{i_r} \right) \right|_{\mu=0} 
  \end{align*}
  coincides with the highest order coefficient of $(b_I)_\mu *
  (\xf_\mu)^I_*$ that, according to the previous lemma and the
  $*$-multiplication formula (\ref{mult}), equals
  \begin{align*}
    \sum_{j=n_{r+1}+1}^{n_r} j!c_{I,j}\Psi_{I,j}.
  \end{align*}
  The equality is therefore proven by the following calculation:
  \begin{align*}
    \left. \frac{\partial^j}{\partial \mu^j} \left(
      (x^1_\mu)^{i_1}(x^2_\mu)^{i_2} \cdots (x_\mu^r)^{i_r} \right)
  \right|_{\mu=0} & = \left. \frac{\partial^j}{\partial \mu^j}
    \prod_{s=1}^r (x^s_0 + x^s_1\mu + \cdots
    x^s_{n_s}\mu^{n_s})^{i_s} \right|_{\mu=0} \\
  & = \left. \frac{\partial^j}{\partial \mu^j}
    \sum_{t=0}^{i_1n_1+\cdots + i_rn_r} \Psi_{I,t}\mu^t
  \right|_{\mu=0}\\  
  & = j!\Psi_{I,j}. 
  \end{align*}
\end{proof}

\subsubsection{The general  complex case}

Since the general solution of (\ref{fmg}) can be decomposed as a sum
(\ref{sumfg}) of solutions to subsystems, one of the main results of
this paper is an immediate consequence of theorem \ref{ps2}:

\begin{Thm3}
  Over the complex numbers, every analytic solution of the matrix
  equation $($\ref{fmg}$)$ is $*$-analytic. 
\end{Thm3}

\subsection{The real case}

The sub-equations of (\ref{fmg}) corresponding to real eigenvalues of
$M$ are treated in the exact same way as in the complex case.
Therefore, it is enough to consider equations corresponding to complex
conjugate pairs of eigenvalues.

We consider only the case when $M$ consists of a single real Jordan
block, i.e., $M$ has the form (\ref{realjordan}). As we have seen in
section \ref{cstmultsec}, there is no restriction to assume that
$\alpha = 0$ and $\beta = 1$. We have already proved that the
corresponding system (\ref{fmg}) can be extended to the parameter
dependent equation (\ref{amuv}). 

In the complex case with one single Jordan block, there was a very
natural choice of ``simple'' solution to construct $*$-power series
form. Namely, $x_\mu = x_0 + x_1\mu + \cdots + x_n\mu^n$ where the
coefficients are given by the coordinate functions. In the current
situation on the other hand, the corresponding function, with
coefficients given by the coordinate functions, is not even a
solution. Instead of searching for an appropriate simple solution,
suitable for constructing power series, we change variables $\xf =
B\fs$, where $B$ is a constant matrix, in such a way that the simple
function
\begin{align}\label{smu}
  s_\mu = s_0 + s_1\mu + \cdots + s_{2n+1}\mu^{2n+1}
\end{align}
becomes a solution. 
\begin{Lem3}\label{lemxbs}
  $s_\mu$ is a solution if $B$ is given by the block structure 
  \begin{align*}
    B_{ij} = \binom{i \!+\! j \!-\! 3}{i - 1}(-\Lambda)^{i+j-2}e_1,
    \quad \begin{array}{l} 
      i = 1, 2, \ldots,  n + 1 \\ 
      j = 1, 2, \ldots, 2(n + 1), 
    \end{array} 
  \end{align*}
  where $B_{ij}$ denotes the $2\times 1$ block in the block-position 
  $(i,j)$ of $B$, and
  \begin{align*}
    \Lambda = \matris{0}{1}{-1}{0},\quad \quad e_1 = \vektor{1}{0}. 
  \end{align*}
\end{Lem3}
\begin{proof}
  By changing variables $\xf = B\fs$, equation (\ref{amuv}) transforms
  into 
  \begin{align*}
    (B^TM^{-1}B^{-T} + \mu I)\nabla V_\mu \equiv 0 \quad \left(mod\; 
      (\mu^2 + 1)^{n+1} \right).  
  \end{align*}
  If we let $X = B^TM^{-1}B^{-T}$ and $(1 + \mu^2)^{n+1} = Z_0 +
  \cdots + Z_{2n+1}\mu^{2n+1} + \mu^{2n+2}$, the equation above can be 
  written as   
  \begin{align*}
    X\nabla V_{i} + \nabla V_{i-1} = Z_i\nabla V_{2n+1}, \quad i = 0,
    1, \ldots, 2n + 1, \quad V_{-1} := 0. 
  \end{align*}
  Thus, if we require $V_\mu = s_\mu$ to be a solution, we see that
  $X$ is uniquely determined as $X = -C^{T}$, where $C$ is the
  companion matrix of the polynomial $(1 + \mu^2)^{n+1}$. In other
  words, we seek a constant matrix $B$ such that $X = -C^{T}$, or
  equivalently, $-M^{-T} = BCB^{-1}$. Thus, from the basic theory of
  rational canonical forms of linear mappings, it is then clear that
  we can choose $B = [\vf\; N\vf\; \cdots\; N^{2n+1}\vf]$, where $N =
  -M^{-T}$ and $\vf$ is any cyclic vector for the matrix $N$.
  Especially, we can choose $v = [1\; 0\; \cdots\; 0]^T$. The proof is
  complete when we show that the powers of $N$ have the following
  block triangular form     
  \begin{align}\label{na}
    (N^a)_{ij} = \left\{ \begin{array}{ll}
          \displaystyle
          \binom{a\!-\!1\!+\!i\!-\!j}{i\!-\!j}(-\Lambda)^{a+i-j} &
          \textrm{if }i \ge j \\  
          & \\
          0 &  \textrm{if }i < j,
        \end{array}\right.
  \end{align}
  where $ i, j = 1, \ldots, n+1$, and $(N^a)_{ij}$ denotes the
  $2\times 2$ block in the block-position $(i,j)$ of $N^a$.  The
  formula (\ref{na}) is true for $a = 0$. Thus, if we assume that it
  is true for some $a \ge 0$, it follows by induction for arbitrary
  $a$:   
  \begin{align*}
    (N^{a+1})_{ij} & = \sum_{k=1}^{n+1} N_{ik}(N^a)_{kj} =
    \sum_{k=j}^{i} (-\Lambda)^{1+i-k}
    \binom{a\!-\!1\!+\!k\!-\!j}{k\!-\!j} (-\Lambda)^{a+k-j} \\  
    & = \sum_{s=0}^{i-j} \binom{a\!+\!s\!-\!1}{s} (-\Lambda)^{a+1+i-j} =
    \binom{a\!+\!i\!-\!j}{i-j}(-\Lambda)^{a+1+i-j}.    
  \end{align*}  
\end{proof}

In the $s$-coordinates, every analytic solution of (\ref{fmg}) can be
expressed as finite sum of $*$-power series 
\begin{align*}
  a_\mu*\sum_{j=0}^\infty c_j(s_\mu)_*^j,
\end{align*}
where $a_\mu$ is a constant solution.

\begin{Thm3}\label{realps}
  Assume that the matrix $M$ consists of a single real Jordan block
  $($\ref{realjordan}$)$, corresponding to a complex conjugate pair of  
  eigenvalues. Then, every analytic solution of $($\ref{fmg}$)$ is
  $*$-analytic. 

  In detail, let $(f, g)$ be an analytic solution of $($\ref{fmg}$)$
  defined by $f + \mathrm{i}g = F_0 + F_1 + \cdots + F_n$, where $F_k$
  is given by $($\ref{Fkreal}$)$ and the corresponding function
  $\phi_k$ has power series representation
  \begin{align*}
    \phi_k(s) = \sum_{j=0}^\infty c_{kj}s^j. 
  \end{align*}
  Then $(f, g)$ can be expressed as a $*$-power series in the
  following way. 
  \begin{align}\label{psreal}
    \sum_{m=0}^n \sum_{j=0}^\infty c_{mj} a_{m,\mu}* (s_\mu)_*^j = f +
    \cdots + g\mu^{2n+1}, 
  \end{align} 
  where $s_\mu$ is given by $($\ref{smu}$)$ and $a_{m,\mu}$ are
  constant solutions given by
  \begin{align*}
    a_{m,\mu} = \sum_{j=0}^m m!(1 + \mu^2)^{n-j} \sum_{k=0}^j
    \sum_{i=0}^m (-1)^k 2^{i-m} \binom{j}{k}b_{ikm} 
  \end{align*}
  where
  \begin{align*}
    b_{ikm} = \left\{
      \begin{array}{ll}\displaystyle
        (-1)^{\frac{m}{2}} \binom{i+2k}{i} & \textrm{if m is even} \\   
        & \\ \displaystyle 
        \mu (-1)^{\frac{m-1}{2}} \binom{i+2k-1}{i} & \textrm{if m is
          odd}
      \end{array}\right.
  \end{align*}
\end{Thm3}

\begin{proof}
  The proof is by induction over $n$. For $n = 0$, $f + \mathrm{i}g =
  \sum_j c_{0j}(x_0 + \mathrm{i}y_0)^j$,
  \begin{align*}
    B = \left[ \binom{-1}{0}e_1 \; \binom{0}{0}(-\Lambda)e_1\right] =
    I. 
  \end{align*}
  Thus, $\xf = \fs$, and since $a_{0,\mu} = b_{000} = 1,$ the $*$-power
  series in the left hand side of (\ref{psreal}) reduces to
  \begin{align*}
    \sum_{j=0}^\infty c_{0j} a_{0,\mu}* (s_\mu)_*^j =
    \sum_{j=0}^\infty c_{0j} (x_0 + y_0\mu)_*^j = f + g\mu.
  \end{align*}
  Hence, the theorem holds for $n = 0$. If $V_{n,\mu}$ denotes the
  $*$-power series in the left hand side of (\ref{psreal}). Then,  
  \begin{align*}
    \sum_{m=0}^n \sum_{j=0}^\infty c_{mj} a_{m,\mu}* (s_\mu)_*^j =
    \sum_{j=0}^\infty c_{nj} a_{n,\mu}* (s_\mu)_*^j + V_{n-1,\mu}(1 +
    \mu^2). 
  \end{align*}
  Thus, if we assume that $V_{n-1,\mu} = \tilde{f} + \cdots +
  \tilde{g}\mu^{2n-1}$ where $\tilde{f} + \mathrm{i}\tilde{g} = F_0 + 
  F_1 + \cdots + F_{n-1}$, it is enough (by induction) to prove that
   \begin{align*} 
     \sum_{j=0}^\infty c_{nj} a_{n,\mu}* (s_\mu)_*^j = f_n + \cdots +
     g_n\mu^{2n+1},  
   \end{align*}
   where $f_n + \mathrm{i}g_n = F_n$. By letting $F_n^{(N)} =
   f_n^{(N)} + \mathrm{i}g_n^{(N)}$ denote the function $F_n$,
   corresponding to $\phi(s) = s^N$, the proof is complete if we can
   prove for any $N$ that 
    \begin{align}\label{amusmuN} 
     a_\mu * (s_\mu)_*^N = f_n^{(N)} + \cdots + g_n^{(N)}\mu^{2n+1},  
   \end{align}   
   where $a_\mu = a_{n,\mu}$. We start by confirming the case when $N = 1$,
   and thereafter we consider the general case $N > 1$. 

   We begin by expressing the functions $f_n^{(1)}$ and $g_n^{(1)}$ in
   the $s$-variables. Let $z_\mu = z_0 + z_1\mu + \cdots z_n\mu^n$,
   then 
   \begin{align*}
     F_n^{(1)} & =  \left. \left( \frac{\partial^n}{\partial
           \mu^n}z_\mu - \sum_{l=1}^n \left( \frac{-\mathrm{i}}{2}
         \right)^l \frac{n!}{(n-l)!}
         \overline{\frac{\partial^{n-l}}{\partial \mu^{n-l}} z_\mu}
       \right) \right|_{\mu=0} \\
    & = \frac{n!}{2^n}\left( 2^nz_n - \sum_{l=1}^n
      2^{n-l}(-\mathrm{i})^l \bar{z}_{n-l} \right).     
   \end{align*}
   Thus, if we let 
   \begin{align*}
     D = \matris{1}{0}{0}{-1},
   \end{align*}
   then, by using lemma \ref{lemxbs} and the identity $\Lambda D = -D
   \Lambda$, we get in matrix notation 
   \begin{align*}
     & \vektor{f_n^{(1)}}{g_n^{(1)}}  = \frac{n!}{2^n}\left(
       2^n\vektor{x_n}{y_n} - \sum_{l=1}^n
       2^{n-l}\Lambda^l\vektor{x_{n-l}}{-y_{n-l}}\right) \\
     & \quad \quad = \frac{n!}{2^n}\left( 2^n\vektor{x_n}{y_n} - \sum_{l=1}^n
       2^{n-l}\Lambda^lD\vektor{x_{n-l}}{y_{n-l}}\right) \\
     &\quad \quad = \frac{n!}{2^n} \left[ -\Lambda^n D \;\; \ldots\;\;
       -2^{n-1} \Lambda D \;\; 2^nI \right] \mathbf{x} \\ 
     & \quad \quad= \frac{n!}{2^n} \left[ -D(-\Lambda)^n \;\; \ldots\;\;
       -2^{n-1}D(-\Lambda) \;\; 2^nI \right] B\mathbf{s} 
\end{align*}
\begin{align*}
     & \quad \quad= \frac{n!}{2^n} \sum_{j=1}^{2n+2} \left( 2^nB_{n+1,j}
       -\sum_{i=1}^n 2^{i-1} D(-\Lambda)^{n-i+1} B_{i,j}
     \right)s_{j-1} \\ 
     & \quad \quad= \frac{n!}{2^n} \!\sum_{j=1}^{2n+2} \!\left( 2^n \binom{n \!+\! j
         \!-\! 2}{n} I  \!-\!\sum_{i=1}^n 2^{i-1} \binom{i \!+\! j
         \!-\!3}{i-1} D  \right) (-\Lambda)^{n+j-1} e_1s_{j-1} \\ 
     & \quad \quad= \sum_{j=0}^{2n+1}\matris{b_j}{0}{0}{d_j}
     (-\Lambda)^{n+j} e_1s_j,  
   \end{align*}
   where
   \begin{align*}
     d_j = \frac{n!}{2^n} \sum_{i=0}^n 2^i \binom{i \!+\! j
       \!-\!1}{i}, \quad \textrm{and } b_j = 2n!\binom{n \!+\!
       j\!-\!1}{n} - d_j.
   \end{align*}
   If we assume that $n$ is even, since $\Lambda^2 = -I$, we obtain 
   \begin{align*}
     & \vektor{f_n^{(1)}}{g_n^{(1)}}  =
     \sum_{j=0}^{2n+1} (-1)^{\frac{n}{2}} \matris{b_j}{0}{0}{d_j}
     (-\Lambda)^{j} e_1s_j \\
     &\quad \quad = \sum_{r=0}^{n} (-1)^{\frac{n}{2}+r} \left( s_{2r}
       \matris{b_{2r}}{0}{0}{d_{2r}} - s_{2r+1}
       \matris{b_{2r+1}}{0}{0}{d_{2r+1}} \Lambda \right) e_1 \\
     &\quad \quad = \sum_{r=0}^{n} (-1)^{\frac{n}{2}+r}
     \vektor{b_{2r}s_{2r}}{d_{2r+1}s_{2r+1}}. 
   \end{align*}
   When $n$ is odd, we get in a similar way
   \begin{align*}
     \vektor{f_n^{(1)}}{g_n^{(1)}} = \sum_{r=0}^{n} (-1)^{\frac{n-1}{2}+r}
     \vektor{-b_{2r+1}s_{2r+1}}{d_{2r}s_{2r}}. 
   \end{align*}

   We now have to confirm that $a_\mu*s_\mu = f_n^{(1)} + \cdots +
   g_n^{(1)}\mu^{2n+1}$. We note that it is enough to prove that the
   coefficients at the highest power of $\mu$ coincide, since the
   other coefficients are then uniquely determined up to irrelevant 
   constant terms. Moreover, we consider only the case when $n$ is
   even, since the case with odd $n$ is completely analogous. Let
   $a_\mu  = \sum_{i=0}^n a_i (1 + \mu^2)^i$, and express $s_\mu$ in
   the following way.  
   \begin{align*}
     s_\mu & = \sum_{r=0}^{2n+1} s_r\mu^r = \sum_{k=0}^n (s_{2k} +
     s_{2k+1}\mu)\mu^{2k} = \sum_{k=0}^n (s_{2k} + s_{2k+1}\mu)(1 +
     \mu^2 - 1)^k \\ 
     & = \sum_{j=0}^n t_j(1 + \mu^2)^j, \quad \textrm{where} \quad
     t_j = \sum_{k=j}^n (-1)^{k-j}\binom{k}{j}(s_{2k} + s_{2k+1}\mu).  
   \end{align*}
   Then, we obtain
   \begin{align*}
     a_\mu * s_\mu & = \left( \sum_{i=0}^n a_i (1 + \mu^2)^i \right) * 
     \left(  \sum_{j=0}^n t_j(1 + \mu^2)^j \right) \\
     & = \sum_{k=0}^n (1 + \mu^2)^k \sum_{r=0}^k a_{k-r}t_r \\ 
     & = \cdots + \mu^{2n+1}\left( \sum_{i=0}^n a_{n-i} \sum_{r=i}^n
       (-1)^{r-i} \binom{r}{i}s_{2r+1}\right) \\
     & = \cdots + \mu^{2n+1}\left( \sum_{r=0}^n s_{2r+1} \sum_{i=0}^r
       (-1)^{r-i} \binom{r}{i} a_{n-i} \right). \\
   \end{align*}
   It is therefore enough, for proving the case $N = 1$, to confirm
   that the following expression vanishes for each $r = 0, 1, \ldots,
   n$. 
   \begin{align*}
     & \sum_{i=0}^r (-1)^i \binom{r}{i} a_{n-i} - (-1)^{\frac{n}{2}}
     d_{2r+1} \\ 
      = & \sum_{i=0}^r (-1)^i \binom{r}{i} n! \sum_{k=0}^i
     \sum_{j=0}^n (-1)^k 2^{j-n} \binom{i}{k} (-1)^{\frac{n}{2}}
     \binom{j\!+\!2k}{j} \\  
     & \quad -  (-1)^{\frac{n}{2}} \frac{n!}{2^n} \sum_{j=0}^n
     \binom{j\!+\!2r}{j}2^j \\
     = & \; (-1)^{\frac{n}{2}} \frac{n!}{2^n} \sum_{j=0}^n 2^j \left(
       \sum_{k=0}^r (-1)^k \binom{j\!+\!2k}{j} \sum_{i=k}^r (-1)^i
       \binom{r}{i} \binom{i}{k} - \binom{j\!+\!2r}{j} \right).    
   \end{align*}
   Thus, the case $N = 1$ is proved by the following calculation.
   \begin{align*}
     & \sum_{i=k}^r (-1)^i \binom{r}{i} \binom{i}{k} =  \sum_{i=k}^r
     (-1)^i \binom{r}{k} \binom{r\!-\!k}{i\!-\!k}  \\ 
      = & \; (-1)^k \binom{r}{k} \sum_{i=0}^{r-k} (-1)^i
     \binom{r\!-\!k}{i}  = 
     \left\{ \begin{array}{ll} 
         (-1)^k &  \textrm{if } r = k \\
         & \\
         0 & \textrm{if } r \ne k.
       \end{array} \right.     
   \end{align*}

   We are now ready to prove (\ref{amusmuN}) for arbitrary $N \ge 1$.
   We have already proven the case when $n = 0$, so it is no
   restriction to assume that $n \ge 1$. Let $\tilde{f}_n^{(N)}$ and
   $\tilde{g}_n^{(N)}$ denote the coefficients at the lowest and
   highest powers of $\mu$ in $a_\mu * (s_\mu)_*^N$, respectively. We
   will prove that  $\tilde{f}_n^{(N)} + \mathrm{i} \tilde{g}_n^{(N)}
   = F_n^{(N)}$ by showing that, when considered as polynomials in
   $x_0$, their coefficients at the power $x_0^{N-1}$ coincide.
   Thereafter we will motivate that this implies that the functions
   must be identical.  

   We consider first the function $F_n^{(N)}$. Since $z_\mu - x_0$ is
   independent of $x_0$, it follows that 
   \begin{align*}
     z_\mu^N & = (z_0 + z_1\mu + \cdots z_n \mu^n)^N  = (x_0 + (z_\mu
     - x_0))^N \\
     & = x_0^N + Nx_0^{N-1}(z_\mu -
     x_0) + \textrm{ (lower order terms in $x_0$ (l.o.t.))} \\
     & = Nx_0^{N-1}z_\mu - (N - 1)x_0^N + \textrm{ (l.o.t.)}
   \end{align*}
   Thus, since $(N - 1)x_0^N$ is independent of $\mu$, we get
   \begin{align*}
     F_n^{(N)} & = \left. \left( \frac{\partial^n}{\partial \mu^n}
         z_\mu^N - \sum_{l=1}^n \left( \frac{-\mathrm{i}}{2} \right)^l
         \frac{n!}{(n-l)!} \frac{\partial^{n-l}}{\partial \mu^{n-l}}
         \bar{z}_\mu^N \right) \right|_{\mu=0} \\
     & =  \left( Nx_0^{N-1} \frac{\partial^n}{\partial \mu^n} z_\mu -
       \sum_{l=1}^{n-1} \left( \frac{-\mathrm{i}}{2} \right)^l
       \frac{n!}{(n-l)!} Nx_0^{N-1} \frac{\partial^{n-l}}{\partial
         \mu^{n-l}} \bar{z}_\mu \right.\\  
     & \quad \quad \left. \left. - \left( \frac{-\mathrm{i}}{2}
         \right)^n n! \bar{z}_\mu^N \right) \right|_{\mu=0} +
     \textrm{ (l.o.t.)} \\
     & = n!x_0^{N-1}\!\left( \!Nz_n \!-\! \sum_{l=1}^{n-1} \left(
         \frac{-\mathrm{i}}{2} \right)^l \! N\bar{z}_{n-l} - \left( 
         \frac{-\mathrm{i}}{2} \right)^n \!(x_0 \!- \!\mathrm{i}Ny_0)\!
     \right) \!+\! \textrm{ (l.o.t.)} \\ 
     & = n!x_0^{N-1}\!\left( Nz_n - \sum_{l=1}^n \left(
         \frac{-\mathrm{i}}{2} \right)^l N\bar{z}_{n-l} - (N - 1)x_0
     \right) + \textrm{ (l.o.t.)} \\ 
     & = x_0^{N-1}\!\left( NF_n^{(1)} - n! (N - 1)x_0 \right) +
     \textrm{ (l.o.t.)}  
   \end{align*}
   On the other hand, from lemma \ref{lemxbs}, it follows that, $s_0 = 
   x_0 + c_0$, where $c_0$ is independent of $x_0$, and that $s_1,
   \ldots s_{2n+1}$ are independent of $x_0$. 
   Thus,
   \begin{align*}
     s_\mu^N = \left( x_0 + (s_\mu - x_0) \right)^N = Nx_0^{N-1}s_\mu
     - (N - 1)x_0^N + \textrm{ (l.o.t.)}, 
   \end{align*} 
   and therefore, using the fact that we already proved the case when
   $N = 1$, we obtain
   \begin{align*}
     a_\mu * (s_\mu)_*^N & = a_\mu * \left( Nx_0^{N-1}s_\mu - (N -
       1)x_0^N + \textrm{ (l.o.t.)} \right) \\ 
     & = x_0^{N-1} \left( Na_\mu*s_\mu - (N - 1)x_0a_\mu \right) +
     \textrm{ (l.o.t.)} \\
     & = x_0^{N-1} \left( Nf_n^{(1)} - (N - 1)a_0x_0 \right.\\
     & \quad \quad \left. + \cdots + \mu^{2n+1} (Ng_n^{(1)} - (N -
       1)a_nx_0)  \right) + \textrm{ (l.o.t.)}  
   \end{align*}
   Hence,
   \begin{align*}
     \tilde{f}_n^{(N)} + \mathrm{i} \tilde{g}_n^{(N)} = x_0^{N-1} \!
     \left( NF_n^{(1)} - (N - 1)(a_0 + \mathrm{i}a_n) x_0 \right) +
     \textrm{ (l.o.t.)} 
   \end{align*}
   we conclude that the coefficients of $\tilde{f}_n^{(N)} +
   \mathrm{i} \tilde{g}_n^{(N)}$ and $F_n^{(N)}$ at the power
   $x_0^{N-1}$ coincide. 

   The last part of the proof is to motivate that $\tilde{f}_n^{(N)} +
   \mathrm{i} \tilde{g}_n^{(N)}$ and $F_n^{(N)}$ must coincide. From
   the $*$-multiplication theorem \ref{realmult}, it follows that
   $(\tilde{f}_n^{(N)}$, $\tilde{g}_n^{(N)})$ must be a solution of
   (\ref{fmg}). Therefore, $\tilde{f}_n^{(N)} + \mathrm{i}
   \tilde{g}_n^{(N)} = F_0 + \cdots + F_n$, where the functions $F_k$
   have the form (\ref{Fkreal}) for some choice of $\phi_k$. Moreover,
   $\tilde{f}_n^{(N)}$ and $\tilde{g}_n^{(N)}$ are homogeneous
   polynomials in the variables $x_0, y_0, \ldots, y_n$. Since $F_k =
   F_k^{(r)}$ will give a polynomial solution which is homogeneous of
   degree $r$, we can therefore conclude that $F_0 = b_0F_0^{(N)}$,
   $F_1 = b_1F_1^{(N)}$, $\cdots$, $F_n = b_nF_n^{(N)}$, for some
   constants $b_0, \ldots, b_n$. Thus, for some polynomial $c$ which
   is constant in $x_0$, we get 
   \begin{align}
     & \tilde{f}_n^{(N)} + \mathrm{i} \tilde{g}_n^{(N)} - F_n^{(N)} =
     b_0F_0^{(N)} + \cdots + b_{n-1}F_{n-1}^{(N)} + (b_n - 1)
     F_n^{(N)} \nonumber \\
     & \quad =  cx_0^N + Nx_0^{N-1} \left( b_0(F_0^{(1)} - x_0) +
       \cdots + b_{n-1}(F_{n-1}^{(1)} - x_0) \right. \nonumber \\ 
     & \quad \quad \quad\left.+ (b_n - 1)(F_n^{(1)} - x_0) \right) +
     \textrm{ (l.o.t.)}   \label{fix0}
   \end{align}
   Finally, the polynomials $F_0^{(1)} - x_0, \ldots, F_n^{(1)} - x_0$
   are linearly independent since, for each $k$, $F_k^{(1)} - x_0$ is
   independent of $y_{k+1}, \ldots, y_n$ with a non-trivial dependence
   on $y_k$.  Therefore, since we know that the
   $x_0^{N-1}$-coefficient in (\ref{fix0}) is zero, we obtain that
   $b_0 = \cdots = b_{n-1} = 0$ and $b_n = 1$. Hence, we conclude
   that  $\tilde{f}_n^{(N)} + \mathrm{i} \tilde{g}_n^{(N)} =
   F_n^{(N)}$, which completes the proof.
\end{proof}

In order to prove that any analytic solution of (\ref{fmg}) is
also $*$-analytic, one would need to cover the last case when the
matrix $M$ consists of several real Jordan blocks corresponding to the
same eigenvalue. This step is technically involved and we have omitted
it here. It is however expected that the result holds in whole
generality, as we have proved in the complex case. We conclude this section by stating a hypothesis for the last remaining case of the equation (\ref{fmg}).

Assume that $M$ is a real matrix with eigenvalues $\pm\mathrm{i}$. It is then no restriction to assume that $M$ has the form $M = \mathrm{diag}(-C^{-T}_1, -C^{-T}_2, \ldots, -C^{-T}_m)$, where each block $C_i$ is the companion matrix for the polynomial $(1 + \mu^2)^{n_i+1}$ with $n_1 \ge n_2 \ge \cdots \ge n_m$. Let $s^1_0, s^1_1, \ldots, s^1_{2n_1+1}, s^2_0, \ldots, s^m_{2n_m+1}$ be the corresponding coordinates. According to theorem \ref{realmult}, the equation (\ref{fmg}) can be extended to the $\mu$-dependent equation (\ref{amuv}). 

\begin{Hyp3}
	The general analytic solution of the equation (\ref{fmg}), with $M = \mathrm{diag}(-C^{-T}_1, -C^{-T}_2, \ldots, -C^{-T}_m)$, is $*$-analytic. In detail, every analytic solution $(f,g)$ of (\ref{fmg}) can be obtained from a $*$-power solution $V_\mu = f + V_1\mu + \cdots + g\mu^{2n+1}$ of (\ref{amuv}) given by
	\begin{align*}
		V_\mu = \sum_{r=1}^m \sum_{I=(i_1,i_2, \ldots , i_r)\in\mathbb{N}^r} \!\!\!\!\!\!(a_I)_\mu * (\mathbf{s}_\mu)^I_*,
	\end{align*}
	where 
	\begin{align}\label{smui}
		(\mathbf{s}_\mu)^I_* := \left( s^1_\mu \right)_*^{i_1} * \left(
      \mu^{n_1-n_2}(s^2_\mu)_*^{i_2} * \cdots * \left(
      \mu^{n_{r-1}-n_r}(s^r_\mu)_*^{i_r} \right)\right). 
	\end{align}
	The $*$-power (\ref{smui}) should be interpreted in a similar way as the corresponding $*$-power (\ref{xI}) in the complex case.
\end{Hyp3}

\section{Conclusions}\label{concsec}

The product of two holomorphic functions is again holomorphic. In
terms of the Cauchy--Riemann equations,
\begin{align*}
  \nabla f = \matris{0}{1}{-1}{0} \nabla g,
\end{align*}
this fact can be expressed through a multiplication (bilinear
operation) $*$ on the solution space
\begin{align*}
  (f, g)*(\tilde{f}, \tilde{g}) = (f\tilde{f} - g\tilde{g}, f\tilde{g} 
  + g\tilde{f}).
\end{align*}
That any holomorphic function is analytic means that any solution
$(f, g)$ of the Cauchy--Riemann equations can be expressed
locally as a convergent power series of a simple (linear in the
variables $x, y$) solution. In a neighborhood of the origin for
example, any solution can be written as
\begin{align*}
    (f, g) = \sum_{r=0}^\infty (a_r, b_r)*(x, y)_*^r, \quad \textrm{where} \quad
    (x, y)_*^r = \underbrace{(x, y)* \cdots *(x, y)}_{r \textrm{ factors}}.
\end{align*}

The main result of this paper is that we establish similar properties
for the more general system
\begin{align}\label{fmg2}
  \nabla f = M \nabla g, 
\end{align}
where $M$ is an arbitrary $n \times n$ matrix with constant entries:
\begin{enumerate}
  \item Any system (\ref{fmg2}) admits a multiplication $*$ on the
    solution space, mapping two solutions to a new solution in a
    bilinear way.
  \item The analytic solutions of (\ref{fmg2}) are characterized by
    the $*$-analytic functions, meaning that every solution of
    (\ref{fmg2}) can be expressed locally through $*$-power series  of
    simple solutions.  
\end{enumerate}

Systems of the form (\ref{fmg2}) constitute only a special case of a quite large family of systems of PDEs which admit $*$-multiplication of solutions \cite{jonasson-2007}. A natural problem, worth studying in the future, is therefore to settle whether the results for the gradient equation (\ref{fmg2}) can be extended to more complex systems of PDEs.

\section*{Acknowledgments}

I would like to thank Prof. Stefan Rauch-Wojciechowski for useful
discussions and comments.

\raggedright
\bibliography{references}
\end{document}